\documentclass[hidelinks,12pt]{amsart}
\renewcommand{\baselinestretch}{1.216}
\usepackage{amsmath}
\usepackage[psamsfonts]{amssymb}
\usepackage{amsfonts}
\usepackage{graphicx}
\usepackage{placeins}

\usepackage{verbatim}
\usepackage{dsfont}
\usepackage{tikz-cd}
\usepackage{mathrsfs}
\usepackage{color}
\usepackage[english]{babel}
\usepackage{indentfirst}
\usepackage{tikz}
\usetikzlibrary{matrix,arrows,decorations.pathmorphing}
\usepackage{algorithm2e}
\usepackage[all]{xy}
\usepackage[T1]{fontenc}
\usepackage{hyperref}
\usepackage{mathtools}
\usepackage{multicol}
\usepackage{enumerate}

\newtheorem{theorem}{Theorem}[section]
\newtheorem{prop}[theorem]{Proposition}
\newtheorem{lemma}[theorem]{Lemma}
\newtheorem{cor}[theorem]{Corollary}

\theoremstyle{definition}
\newtheorem{example}[theorem]{Example}
\theoremstyle{plain}

\newtheorem{remark}[theorem]{Remark}

\renewcommand{\hat}[1]{\widehat{#1}}
\numberwithin{equation}{section}

\DeclareMathOperator{\R}{\mathbb{R}}

\def\ep{\epsilon}

\def\R{\mathbb{R}}
\def\Z{\mathbb{Z}}
\def\N{\mathbb{N}}

\def\C{\mathbb{C}}

\def\T{\mathbb T}

\def\cal R{\mathcal R}

\newcommand{\Sh}{\operatorname{Sh}}
\newcommand{\Vol}{\operatorname{Vol}}
\newcommand{\Ham}{\operatorname{Ham}}
\newcommand{\std}{\mathrm{std}}
\newcommand{\dBM}{d_{\mathrm{BM}}}
\newcommand{\Wred}{\mathcal W}
\newcommand{\Ddiag}{\mathcal D}
\newcommand{\Cred}{\mathcal C}
\newcommand{\Tpm}{\mathcal T_{\pm}}
\newcommand{\Tplus}{\mathcal T_{+}}
\newcommand{\Tminus}{\mathcal T_{-}}
\newtheorem{mainthm}{Theorem}

\newcommand{\into}{\hookrightarrow}

\newcommand{\cECH}[1]{c_{#1}^{\mathrm{ECH}}}

\allowdisplaybreaks[1]
\usepackage[letterpaper,textwidth=6.68in,textheight=8.5in,centering]{geometry}

\author{Richard Hind}
\email{hind.1@nd.edu}
\address{Department of Mathematics, University of Notre Dame, Notre Dame, IN 46556, USA}

\author{Jun Zhang}
\email{jzhang4518@ustc.edu.cn}
\address{The Institute of Geometry and Physics, University of Science and Technology of China, 96 Jinzhai Road, Hefei Anhui, 230026, China}
 
\title{The shape invariant and Banach--Mazur geometry}

\begin{document}

\begin{abstract}
We prove an intersection theorem for Lagrangian tori that extends a recent result of Hind and Kerman. Together with a Hamiltonian folding construction, this yields computations of the reduced shape invariant for a class of four-dimensional toric domains, extending earlier work of Hind and Zhang to non-integer ellipsoids with axis ratio greater than two. We also introduce the polydisk shape invariant and obtain partial computations by adapting the proof of our intersection theorem. Applications include sharp polydisk-to-ellipsoid embedding criteria in new parameter ranges and obstructions not detected by ECH capacities. We also determine when the shape invariant commutes with finite intersections of graphical toric domains. Moreover, we construct a convex domain at positive Banach--Mazur distance from all toric convex and concave domains. Finally, we prove that the space of star-shaped domains is not locally compact, with respect to the topology induced by the Banach--Mazur distance. 
\end{abstract}

\maketitle

\begingroup
\renewcommand{\baselinestretch}{1}

\tableofcontents
\endgroup

\section{Introduction}

Quantitative symplectic topology studies how the geometry of a domain constrains its symplectic embeddings. Gromov's nonsqueezing theorem \cite{Gromov85} shows that these constraints do not come from volume alone. Various symplectic capacities were invented to provide further obstructions. Most symplectic capacities are either a single number or a sequence of scalars, and they encode information about a domain from different perspectives. The shape invariant, originating from the work of Sikorav and Eliashberg \cite{Sikorav89,Sikorav91,Eliashberg91}, implements this idea through the area classes of embedded Lagrangian tori in a given domain, and it results in multi-dimensional regions instead of scalars. 

The purpose of this paper is to develop the reduced shape invariant in dimension four and to use it in embedding problems and Banach--Mazur geometry (in the symplectic setting, introduced by Ostrover--Polterovich around 2015). Our starting point is a new intersection theorem for Lagrangian tori, extending a recent result of Hind--Kerman \cite{HK26} beyond Hamiltonian images of product tori. Combined with the Hamiltonian folding constructions of \cite{HZ23}, it computes the reduced shape invariant of a wide class of toric domains. In particular, it removes the integrality restriction from Hind--Zhang's ellipsoid computations of \cite{HZ26} when the ratio of the axes is at least two. We also consider the \emph{polydisk shape invariant}. A further application of the method underlying the intersection theorem determines a portion of this invariant which leads to obstructions which do not follow from the reduced shape invariant. In particular we get sharp obstructions for polydisk embeddings into ellipsoids in new parameter ranges, considerably generalizing the so-far-all-known obstructions (say, from Hutchings \cite{Hutchings16} or from Digiosia--Nelson--Ning--Weiler--Yang \cite{DNNWY22}). We also study the shape under finite intersections: our new intersection theorem gives an exact commutation criterion for graphical toric domains and computations for aligned ellipsoids. Finally, we apply shape invariant information to distinguish convex domains from convex and concave toric domains quantitatively, and obtain an example that is comparable with the counterexample from Haim-Kislev--Ostrover in \cite{HKO2026}. Moreover, we investigate the local geometry of the space of star-shaped domains, with respect to the Banach--Mazur geometry. 

\subsection{Notation and shape invariants}

We work in $\R^4=\C^2$ with $\omega_\std\coloneqq\sum_{j=1}^2 dx_j\wedge dy_j, \lambda_\std\coloneqq\frac12\sum_{j=1}^2(x_j\,dy_j-y_j\,dx_j)$ and $z_j\coloneqq x_j+iy_j$. The moment map is
\[
 \mu\colon\C^2\longrightarrow\R^2_{\geq0},
 \qquad \mu(z_1,z_2)\coloneqq(\pi|z_1|^2,\pi|z_2|^2).
\]
For a region $\Omega\subset\R^2_{\geq0}$, write $X_\Omega\coloneqq\mu^{-1}(\Omega)$. For a continuous function $f\colon[0,a]\to[0,\infty)$, with $a>0$ and $f(u)>0$ for $u<a$, denote the graphical moment region by
\[
 \Omega_f\coloneqq\{(u,v)\in\R^2_{\geq0}\mid u<a,\ v<f(u)\}.
\]
For $0<a<\infty$ and $a\leq b\leq\infty$, a four-dimensional symplectic ellipsoid is denoted by 
\[
 E(a,b)\coloneqq\left\{z\in\C^2\;\middle|\;
 \frac{\pi|z_1|^2}{a}+\frac{\pi|z_2|^2}{b}<1\right\},
\]
where $1/\infty\coloneqq0$. In particular, the symplectic ball and symplectic cylinder are the special cases: $B^4(a)\coloneqq E(a,a)$ and $Z^4(a)\coloneqq E(a,\infty).$ For $0<c\leq d<\infty$, the symplectic polydisk is denoted by 
\[
 P(c,d)\coloneqq\left\{z\in\C^2\;\middle|\;
 \pi|z_1|^2<c,\ \pi|z_2|^2<d\right\}.
\]
The product torus with action pair $(r,s)$ is $L(r,s)\coloneqq\mu^{-1}(r,s)$, where $r,s>0$. Throughout, the axes of ellipsoids and polydisks are ordered as above. For finite $b$, an ellipsoid is called \emph{integral} when $b/a\in\N_{\geq2}$, following \cite{HZ23,HZ26}; a \emph{non-integer ellipsoid} is one whose axis ratio is not an integer.

Except for the explicitly indicated cylinder, domains are open and bounded, and star-shaped domains are star-shaped with respect to the origin. We allow the standard piecewise smooth toric models; statements involving their boundaries can equivalently be formulated using arbitrarily small smooth enlargements. An arrow $X\into Y$ denotes a symplectic embedding. A Hamiltonian embedding of a torus (or an embedded torus) means an image under a Hamiltonian diffeomorphism of the ambient $(\C^2, \omega_{\rm std})$. For $t>0$ we use the following convention of rescaling:
\begin{equation}\label{eq:scaling}
 tX\coloneqq\{\sqrt{t}\,z\mid z\in X\}; \qquad \Vol(tX)=t^2\Vol(X) \,\,\mbox{where}\,\, \Vol(X)\coloneqq\frac12\int_X\omega_\std^2.
\end{equation}
Thus $tE(a,b)=E(ta,tb)$ and $tP(c,d)=P(tc,td)$.

\medskip

Let us recall the definition of the reduced shape invariant (introduced in \cite{HZ26}). Let $L\subset X$ be an embedded Lagrangian torus. Its area class and Maslov class are the homomorphisms
\[
 \mathcal A_L\colon H_1(L;\Z)\longrightarrow\R \,\,\,\,\mbox{and}\,\,\,\,
 m_L\colon H_1(L;\Z)\longrightarrow2\Z.
\]
The first records the symplectic areas of disks bounded by loops on $L$ (which can be expressed as $\mathcal A_L(e)\coloneqq\int_e\lambda_\std$); the second records the winding of the tangent Lagrangian planes. A Maslov-two integral basis $(e_1,e_2)$ need not be unique. For example, replacing it by $(2e_1-e_2,e_1)$ preserves the Maslov indices and changes the ordered action pair from $(r,s)$ to $(2r-s,r)$. Consequently, the area pair of a parametrized torus contains a certain ambiguity that should be removed before comparing domains.

The normalization in \cite[Lemma~1.1]{HZ26}, based on the existence results of holomorphic disks from Viterbo and Polterovich \cite{Viterbo90,Polterovich91}, resolves this ambiguity. Every embedded Lagrangian torus in $\R^4$ admits a Maslov-two integral basis with positive ordered actions. If the torus is nonmonotone, there is a unique such basis satisfying
\[
 0<2\mathcal A_L(e_1)\leq\mathcal A_L(e_2).
\]
If the torus is monotone, the two actions are equal and their common value is independent of the choice of Maslov-two basis. For brevity, set
\[
 \Wred\coloneqq\{(r,s)\mid 0<2r\leq s\},\qquad
 \Ddiag\coloneqq\{(r,r)\mid r>0\},\qquad
 \Cred\coloneqq\Wred\cup\Ddiag.
\]
The \emph{reduced shape invariant} is (see (4) in \cite{HZ26}): 
\begin{equation*}
 \Sh^+(X)\coloneqq\left\{\bigl(\mathcal A_L(e_1),\mathcal A_L(e_2)\bigr)\;\middle|\;
 \begin{array}{l}
 L\subset X\text{ an embedded Lagrangian torus}\\
 (e_1,e_2)\text{normalized in \cite[Lemma~1.1]{HZ26}}
 \end{array}\right\}\subset\Cred.
\end{equation*}
This shape invariant admits a variant, called the \emph{Hamiltonian shape invariant}, defined by 
\begin{equation*}
 \Sh_H^+(X)\coloneqq\left\{(r,s)\;\middle|\;
 \begin{array}{l}
 0<r\leq s,\ \Phi(L(r,s))\subset X\\
 \text{for some }\Phi\in\Ham(\C^2)
 \end{array}\right\}.
\end{equation*}
Here, the order $r\leq s$ only removes the symmetry interchanging the two factors, and it does not impose the stronger region constraint $s\geq2r$. By definition, on $\Cred$ there is an inclusion
\begin{equation*}
 \Sh_H^+(X)\cap\Cred\subset\Sh^+(X).
\end{equation*}

Note that, for arbitrary domains in $\R^4$, this inclusion could be strict. 

\begin{example}[Maslov-two base change]\label{ex:reduced-actions}
The product torus $L(2,3)$ has reduced area pair $(1,2)$: if $(e_1,e_2)$ is its standard basis, the Maslov-two basis $(2e_1-e_2,e_1)$ has actions $(1,2)$. However, $L(2,3)$ is not Hamiltonian isotopic to $L(1,2)$; see \cite{Chekanov96} and \cite[Remark~1.7]{HZ26}. Thus the reduced area pairs do not determine the Hamiltonian isotopy class, and the reduced shape invariant does not automatically determine the full Hamiltonian shape invariant. As a matter of fact, one can show that when $X$ is a sufficiently small neighborhood of $L(2,3)$, then $(1,2) \notin  \Sh_H^+(X)$ while $(1,2) \in \Sh^+(X)$ by the base change as mentioned above. 
\end{example}

Both invariants are monotone under symplectic embeddings of star-shaped domains and rescaling under \eqref{eq:scaling}; see \cite[Proposition~6.1]{HZ26} and \cite{HZ23}. In other words,
\[
 X\into Y\ \Longrightarrow\ \Sh^+(X)\subset\Sh^+(Y),
 \qquad \Sh^+(tX)=t\Sh^+(X),
\]
with analogous statements for $\Sh_H^+$. With $\Omega^*\coloneqq\{(s,r)\mid (r,s)\in\Omega\}$, the fibers (as product Lagrangian tori) yield
\[
 (\Omega\cup\Omega^*)\cap\{0<r\leq s\}\subset\Sh_H^+(X_\Omega).
\]
The main task is to understand how far the shape invariant extends beyond the moment image.

\begin{example}[Hamiltonian shape of standard domains]\label{ex:standard-shapes}
The computations of Hind--Opshtein \cite{HO20} and Hind--Zhang \cite{HZ26} exhibit both bounded moment regions and unbounded strips. For $R>0$, the ball satisfies
\[
 \Sh_H^+(B^4(R))
 =\bigl(\{(r,s)\mid r+s<R\}\cup\{(r,s)\mid r<R/3\}\bigr)
 \cap\{(r,s)\mid 0<r\leq s\}.
\]
For $0<a\leq b<\infty$ with $b/a\in\N_{\geq2}$, \cite[Theorems~1.5 and~1.8]{HZ26} gives
\[
 \Sh_H^+(E(a,b))=
 \left(\left\{(r,s)\;\middle|\;\frac r a+\frac s b<1\right\}
 \cup\left\{0< r<\frac{a}{2}\right\}\right)\cap\{(r,s)\mid 0<r\leq s\}.
\]
Both domains are ellipsoids, since $B^4(R)=E(R,R)$, but their unbounded strips have different widths: $R/3$ for the ball and $a/2$ for the integral ellipsoid with $b\geq2a$. 

For $0<c\leq d<\infty$, the polydisk satisfies
\[
 \Sh_H^+(P(c,d))
 =\bigl(\{(r,s)\mid r<c,\ s<d\}\cup\{(r,s)\mid r<c/2\}\bigr)
 \cap\{(r,s)\mid 0<r\leq s\}.
\]
In each case above, $\Sh^+$ is obtained by intersecting the right-hand side with $\Cred$. For illustrations, together with the symplectic cylinder case, see the diagram in \cite[p.~6]{HZ26}. The unbounded strip records the possibility of placing (via Hamiltonian isotopies) product tori with arbitrarily large second action in a bounded domain when their first action is sufficiently small. Beyond that strip, the moment boundary serves as an obstruction.
\end{example}

\subsection{Rigidity and flexibility}

The recent work of Hind--Kerman \cite{HK26} provides a more fundamental embedding obstruction: if $a/2\leq r<a\leq s$, 
and $\Phi\in\Ham(\C^2)$ satisfies $\Phi(L(r,s))\subset Z^4(a)$, then
\begin{equation*}
 \Phi(L(r,s))\cap \left(\bigcup_{t\geq s} L(r,t)\right) \neq\varnothing.
\end{equation*}
This is \cite[Corollary~1.8]{HK26}, obtained from their intersection theorem in a polydisk. Unlike a capacity inequality, it locates an intersection with a specified family of fibers. Our first main result extends this intersection statement, as a rigidity result, to arbitrary Lagrangian tori in the non-monotone region $\Wred$.

\begin{mainthm}\label{thm:A}
Let $a>0$ and let $\varphi\colon\T^2\into Z^4(a)$ be a Lagrangian embedding. Suppose that $H_1(\varphi(\T^2);\Z)$ has a Maslov-two integral basis with area pair $(r,s)$ satisfying $a \leq 2r \leq s$. Then
\[
 \varphi(\T^2)\cap \left(\bigcup_{t\geq s}L(r,t)\right)\neq\varnothing.
\]
\end{mainthm}

The novelty is that the embedded Lagrangian torus $\varphi(\T^2)$ need not be Hamiltonian isotopic to a product torus. This is essential for an application to $\Sh^+$ rather than only to $\Sh_H^+$. Note that this result is not a formal consequence of \cite{HK26}: their argument uses Hamiltonian-isotopy information, and \cite[Remark~2.11]{HK26} explains why area and Maslov data cannot replace it without additional restrictions. Our hypothesis $s\geq2r$ is precisely the normalization relevant to the non-monotone part of the reduced shape invariant. Moreover, the restriction $s\geq2r$ has a useful arithmetic interpretation. Every Maslov-two class has the form $(1-k)e_1+ke_2$. So it has area $r+k(s-r)$ for some $k\in\Z$. When $s\geq2r$, every positive area of this form is at least $r$, which provides a key point in the argument. Nevertheless, establishing the desired conclusion in Theorem~\ref{thm:A} still requires holomorphic curve constructions and neck-stretching arguments as in \cite{HK26}. The proof of Theorem~\ref{thm:A} will be provided in Section~\ref{sec:proof-A}.

\medskip

For the opposite, flexibility, we use the Hamiltonian folding construction of \cite[Theorem~5.2(2)]{HZ23}, in the following form. Its short proof is given in Section~\ref{sec:proof-B}.

\begin{mainthm}\label{thm:B}
Let $\rho>0$, and set
\[
 T_\rho\coloneqq\operatorname{conv}\{(0,0),(0,2\rho),(\rho,2\rho),(2\rho,0)\} \,\,\,\,\mbox{and}\,\,\,\, Q_\rho\coloneqq\mu^{-1}(T_\rho).
\]
Then for every $0<r<\rho$ and every $s\geq r$, a Hamiltonian image of $L(r,s)$ is contained in the interior of $Q_\rho$. In particular, if an open domain $X\subset\C^2$ contains a neighborhood of $Q_\rho$, then
\[
 \{(r,s)\mid 0<r \leq \rho,\ r\leq s\}\subset\Sh_H^+(X).
\]
\end{mainthm}

The construction in \cite{HZ23} controls an entire Hamiltonian isotopy, and was developed to lift paths in the shape invariant and to detect Hamiltonian knottedness. Here, we only use its endpoint to realize an infinite strip of area pairs. Note that, from the flexibility side, Theorem~\ref{thm:B} explains the different strip widths in Example~\ref{ex:standard-shapes}. The quadrilateral $T_\rho$ lies in the moment image of $B^4(R)=E(R,R)$ when $3\rho<R$, while it lies in the moment image of $E(a,b)$ with $b\geq2a$ when $2\rho<a$. Similarly, $T_\rho$ lies in the moment rectangle of $P(c,d)$ when $2\rho<c\leq d$, giving width $c/2$. Corollary \ref{cor:ellipsoid-shape} below computes the reduced shape invariant of general (including non-integer) ellipsoids $E(a,b)$, where the width of the unbounded strip is still $\frac{a}{2}$ whenever $\frac{b}{a} \geq 2$.  

Curiously, our methods for computing both the reduced shape invariant and the Hamiltonian shape invariant of $E(a,b)$ fail when $1< \frac{b}{a}< 2$. To simplify the notation, up to rescaling, let us consider $E(1,x)$ where $1<x<2$. In this case, there is still an unbounded strip with width denoted by $w(x)$, precisely
\[
w(x) = \sup \{ w \, | \, L(r,s) \hookrightarrow E(1,x) \, \, \mathrm{for} \, \, \mathrm{all} \, \, r \le w, \, s \ge r \}.
\]

It is not clear at all what $w(x)$ looks like (although it is nondecreasing, and we know $\lim_{x \to 1^+} w(x) = \frac{1}{3}$ and $\lim_{x \to 2^-} w(x) = \frac{1}{2}$). However, Theorem \ref{thm:B} helps to narrow down the width $w(x)$ in a certain region. 

\begin{cor} \label{cor-strip-width}
For the symplectic ellipsoid $E(1,x)$ where $1<x<2$, we have  
\[ \frac{x}{x+2} \leq w(x) \leq \min\left\{\frac{1}{2}, \frac{x}{3} \right\} \]
where $w(x)$ denotes the width of the unbounded strip in $\Sh_H^+(E(1,x))$.
\end{cor}

The proof of Corollary \ref{cor-strip-width} will also be given in Section~\ref{sec:proof-B}. Any further improvement of the estimation of $w(x)$ requires either a refined argument of pseudo-holomorphic curves (to reduce the upper bound in Corollary \ref{cor-strip-width}) or a new folding construction (to bump up the lower bound in Corollary \ref{cor-strip-width}). 

\subsection{Reduced shape and non-integer ellipsoids}

Combining Theorems~\ref{thm:A} and~\ref{thm:B}, one can compute the reduced shape invariant of toric domains whose moment regions lie below sufficiently tall concave graphs. The proof of the following corollary, which lies in Section~\ref{sec:proof-B}, makes explicit the roles of the intersection theorem and the Hamiltonian construction.
\begin{cor}\label{cor:graph-shape}
Let $a>0$ and let $f\colon[0,a]\to[0,\infty)$ be a continuous,
nonincreasing, concave function such that $f(u)>0$ for $0\leq u<a$
and $f(a/2)>a$. Then
 \begin{equation}\label{eq:graph-shape}
\Sh^+(X_{\Omega_f}) =\Sh_H^+(X_{\Omega_f})\cap\Cred = \left(\Omega_f\cup \{(r,s)\mid 0<r<a/2,\ s>0\}\right)\cap\Cred.
\end{equation}
\end{cor}

\begin{figure}[!htbp]
  \centering
  \includegraphics[scale=0.8]{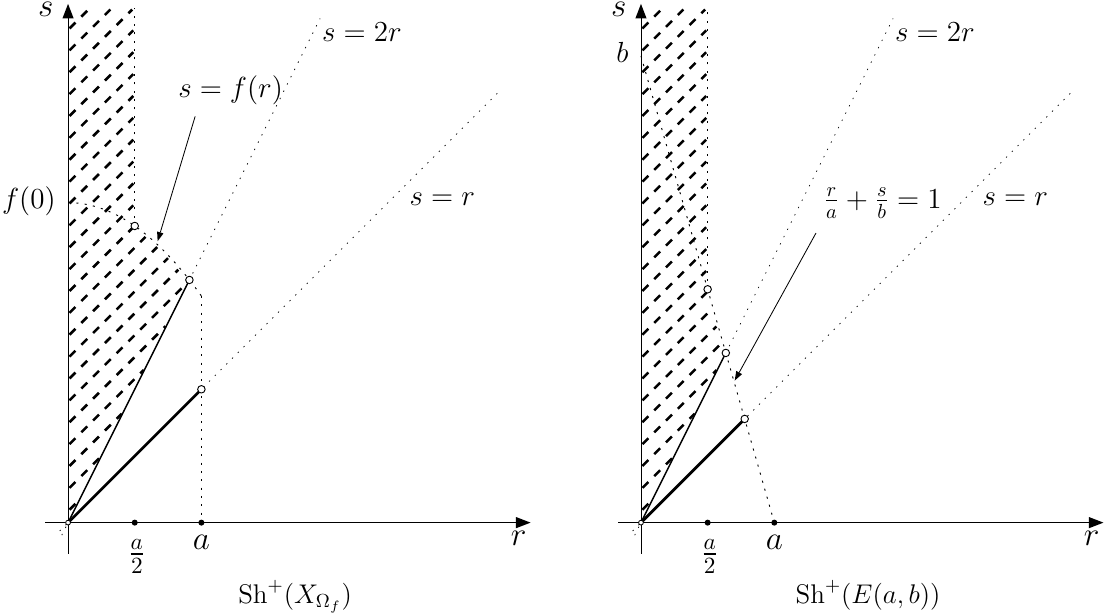}
  \caption{The reduced shape invariants in Corollaries~\ref{cor:graph-shape}
  and~\ref{cor:ellipsoid-shape}.}
  \label{fig:graph-ellipsoid-shapes}
\end{figure}

For an illustration of a function in the hypothesis of Corollary \ref{cor:graph-shape}, see the left picture in Figure \ref{fig:graph-ellipsoid-shapes}. For an ellipsoid, $f(u)=b(1-u/a)$ and the hypothesis $f(a/2)>a$ is equivalent to $b>2a$. Combining this with the previously known case $b=2a$ in \cite{HZ26} gives the following result.

\begin{cor}\label{cor:ellipsoid-shape}
Let $0<a\leq b<\infty$ with $b\geq2a$. Then
\begin{equation}\label{eq:new-ellipsoid-shape}
 \Sh^+(E(a,b))=
 \left(\left\{(r,s)\;\middle|\; \frac r a+\frac s b<1\right\}
 \cup\left\{0<r<\frac{a}{2}\right\}\right)\cap\Cred.
\end{equation}
\end{cor}

Again, for an illustration, see the right picture in Figure \ref{fig:graph-ellipsoid-shapes}. Note that Corollary \ref{cor:ellipsoid-shape} covers non-integer ellipsoids and answers a question explicitly raised in \cite[Remark~1.9]{HZ26} (for the reduced shape invariant). It does not, by itself, compute the Hamiltonian shape invariant on the additional region $r<s<2r$. This distinction is also relevant to the discussion of non-integer ellipsoids in \cite[Corollary~1.20]{HK26}.

\subsection{Polydisk shape}

The reduced shape invariant measures individual tori through normalized area pairs. In an embedding problem, however, a torus may occur as part of an embedded polydisk, and this additional structure can impose stronger constraints. To initiate the study in this direction, define
\begin{equation}\label{eq:polydisk-shape}
 \Sh^P(X)\coloneqq\left\{(r,s)\, \bigg| \, 0<r\leq s,\ \overline{P(r,s)}\into X\right\},
\end{equation}
where an embedding of the closed polydisk is understood to be defined on a neighborhood. 
The polydisk shape is investigated by the two authors together with Ahn and Kerman in \cite{AKHZ}.

Like other versions of shape, the polydisk shape invariant is monotone and conformal. By restricting an embedding to its corner torus and using the ambient extension for embeddings of a star-shaped domain, we have
\begin{equation*}
 \Sh^P(X)\subset\Sh_H^+(X),
 \qquad
 \Sh^P(X)\cap\Cred\subset\Sh_H^+(X)\cap\Cred\subset\Sh^+(X).
\end{equation*}
The inclusion $\Sh^P(X)\subset\Sh_H^+(X)$ can be strict. For a bounded domain, volume gives
\begin{equation*}
 (r,s)\in\Sh^P(X)\ \Longrightarrow\ rs\leq\Vol(X).
\end{equation*}
Thus $\Sh^P(X)$ cannot contain any vertical strip of positive width extending to infinity (recall that the Hamiltonian shape invariant of a polydisk does contain such a strip, as in Example \ref{ex:standard-shapes}).

\begin{example}[Polydisk shape of a cylinder]\label{ex:cylinder-polydisk}
Let $a>0$. Nonsqueezing and the standard inclusion give
\[
 \Sh^P(Z^4(a))=\{(r,s)\mid 0<r<a,\ r\leq s\}.
\]
On the other hand, for any given bounded polydisk, Schlenk's folding constructions, see \cite{Schlenk05}, can enlarge $\Sh^P$ beyond the moment rectangle, but the volume bound continues to constrain this enlargement. 
\end{example}

Recall that a toric domain (as well as its moment image) is called \emph{monotone} when the profile curve of its moment region is the graph of a decreasing function. 
The work \cite{AKHZ} proves the following partial computation of the polydisk shape invariant. 
\begin{theorem}[see \cite{AKHZ}] \label{cor:polydisk-rigidity}
Let $a>0$ and $X_\Omega\coloneqq\mu^{-1}(\Omega)$ be a monotone toric domain contained in $Z^4(a)$. Then, writing $\Omega^*\coloneqq\{(s,r)\mid(r,s)\in\Omega\}$, we have
\begin{equation}\label{eq:polydisk-rigidity}
 \Sh^P(X_\Omega)\cap\{r\geq a/2\}
 = (\Omega\cup\Omega^*)\cap\{0<r\leq s,\ r\geq a/2\}.
\end{equation}
\end{theorem}
For an illustration of the region on the righthand side of (\ref{eq:polydisk-rigidity}), see the following Figure \ref{fig:polydisk-shape-reflection}. Note that the union with the reflection $\Omega^*$ could contain more points than the original $\Omega$. 
\begin{figure}[!htbp]
  \centering
  \includegraphics[scale=0.64]{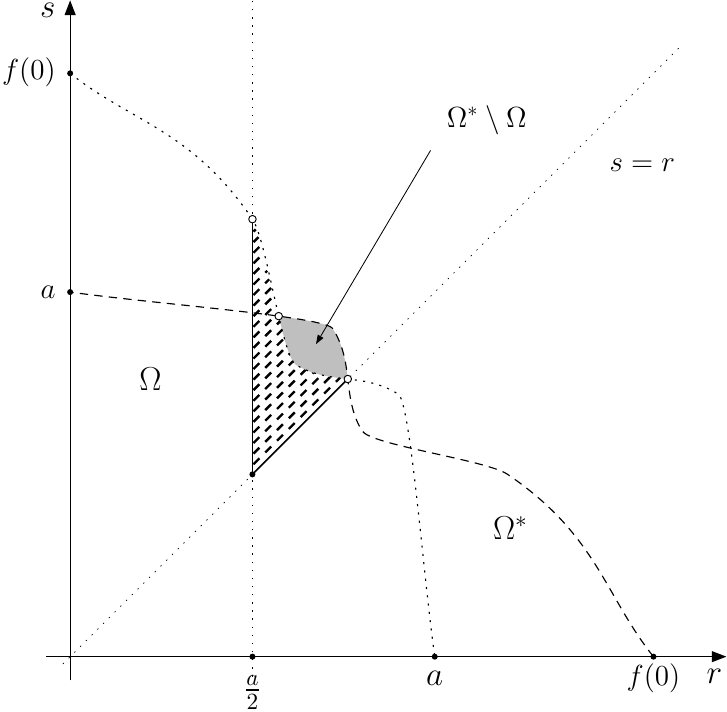}
  \caption{The region $(\Omega\cup\Omega^*)\cap\{r\geq a/2,\ r\leq s\}$ in
  Corollary~\ref{cor:polydisk-rigidity}.}
  \label{fig:polydisk-shape-reflection}
\end{figure}


The argument below proves this rigidity when $s\geq a$, that is, in the case when obstructions can be derived from studying embeddings of the singular torus $L(r,s) \subset \partial P(r,s)$. 

\begin{proof}[Proof of Theorem \ref{cor:polydisk-rigidity}, in the range $s\geq a$]
Suppose that $(r,s)\in\Sh^P(X_\Omega)$ with $r\geq a/2$ and $s\geq a$.
Since the embedding is defined on a neighborhood of $\overline{P(r,s)}$,
nonsqueezing gives $r<a$. Its restriction to the corner torus $L(r,s)$ extends to a
Hamiltonian diffeomorphism of $\C^2$. Hence \cite[Corollary~1.8]{HK26}
implies that its image meets $L(r,t)$ for some $t\geq s$.
Thus $(r,t)\in\Omega$, and monotonicity of the profile gives $(r,s)\in\Omega$.
Conversely, for $(r,s)\in\Omega\cup\Omega^*$ the closed polydisk embeds by
inclusion, with the coordinate swap when needed; openness supplies a
neighborhood. This proves \eqref{eq:polydisk-rigidity} after intersecting
both sides with $\{s\geq a\}$.
\end{proof}

\begin{example}[Polydisks in an ellipsoid]\label{ex:ellipsoid-polydisk}
For every finite $x\geq1$, Theorem~\ref{cor:polydisk-rigidity} yields
\begin{equation}\label{eq:polydisk-ellipsoid-region}
 \Sh^P(E(1,x))\cap\{r\geq1/2\}
 =\left\{(r,s)\;\middle|\; \frac12\leq r\leq s,\ r+\frac s x<1\right\}.
\end{equation}
This statement concerns all axis ratios in the stated range, not merely integer ratios. The geometric restriction is on the size of the smaller polydisk factor relative to the first ellipsoid axis.
\end{example}

\begin{remark} For some cases, say $B^4(R)$, by comparing the Hamiltonian shape invariant of $P(1,a)$ and $B^4(R)$ as in Example \ref{ex:standard-shapes} (back to \cite{HO20}) or by \cite[Theorem~1.5]{Hutchings16}, we know more information than the one given in Theorem \ref{cor:polydisk-rigidity}. Explicitly, we have 
\[
 \Sh^P(B^4(R))\cap\{r\geq R/3\}
 = \mu(B^4(R)) \cap\{0<r\leq s,\ r\geq R/3\}
\]
improving the range from $r\geq R/2$ (from (\ref{eq:polydisk-ellipsoid-region})) to $r\geq R/3$. Indeed, the Hamiltonian shape invariant comparison implies that $P(r,s) \hookrightarrow B^4(R)$ (which is equivalent to $P(1, \frac{s}{r}) \hookrightarrow B^4(\frac{R}{r})$), when restricted to the region $\{r \leq s \leq 2r\}$, if and only if $r+s \leq R$. Moreover, on the boundary line $r+s=R$, the condition $\{r\leq s \leq 2r\}$ precisely corresponds to $\{\frac{R}{3} \leq r \leq \frac{R}{2}\}$. For the ball, the region in (\ref{eq:polydisk-rigidity}) contains no point with $r=R/2$; the point $(R/2,R/2)$ lies only on the boundary $r+s=R$. 
\end{remark}

\subsection{Sharp embedding criteria}

Write
\[
 c_\beta(a)\coloneqq\inf\{c>0\mid P(1,a)\into E(c,\beta c)\},
 \qquad a,\beta\geq1.
\]
The inclusion of moment images always gives $c_\beta(a)\leq1+\frac a\beta$. Interestingly, rescaling \eqref{eq:polydisk-ellipsoid-region} suggests that this also provides the obstruction, i.e., 
\begin{equation} \label{eq:linear-obstruction}
 a+\beta\leq\beta c.
\end{equation}
Indeed, we can rephrase \eqref{eq:polydisk-ellipsoid-region} in the following way: suppose $\overline{P(r,s)} \hookrightarrow E(1,x)$ (which is equivalent to $\overline{P(1, \frac{s}{r})} \hookrightarrow E(\frac{1}{r}, \frac{x}{r})$) and set $a\coloneqq \frac{s}{r}$, $c\coloneqq \frac{1}{r}$, and $\beta\coloneqq x$. Then when $c\leq 2$ (hence, $r \geq \frac{1}{2}$), we must have 
\[ a+\beta < \beta c \,\,\,\, \mbox{or equivalently $r + \frac{s}{x} < 1$.} \]
We also have (\ref{eq:linear-obstruction}) for certain $c>2$, subject to Theorem~\ref{cor:polydisk-rigidity}. This leads to the following embedding obstruction result. 

\begin{mainthm}\label{thm:C}
Let $a,\beta,c\geq1$. The equivalence
\begin{equation}\label{eq:C}
 P(1,a)\into E(c,\beta c)
 \quad\Longleftrightarrow\quad
 a+\beta\leq\beta c
\end{equation}
holds whenever at least one of the following conditions is satisfied:
\begin{enumerate}[(i)]
\item $1\leq c<2$;
\item $\beta\geq2$ and $1\leq a\leq\beta$;
\item $\beta \in \R \backslash (1,2)$ and $1\leq a\leq2$.
\end{enumerate}
\end{mainthm}

\begin{figure}[!htbp]
  \centering
  \includegraphics[width=.78\textwidth]{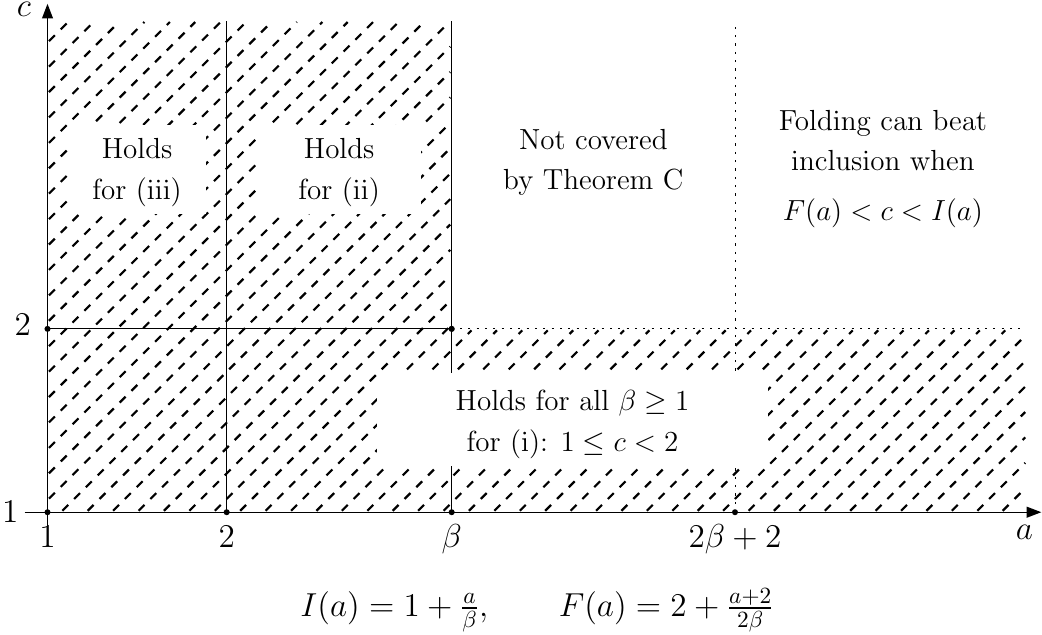}
  \caption{Parameter ranges in Theorem~\ref{thm:C}.
  For $a>2\beta+2$, folding improves the inclusion bound;
  the two bounds $F(a)$ and $I(a)$ are displayed.}
  \label{fig:theorem-c-parameter-ranges}
\end{figure}

The proof is discussed in Section~\ref{sec:proof-C} and Figure \ref{fig:theorem-c-parameter-ranges} lists the involving parameter ranges. As a comparison, \cite[Theorem~1.11(2)]{HZ26} proves \eqref{eq:C} for integral $\beta\geq2$ and $1\leq c\leq2$. Hutchings \cite[Theorem~1.5]{Hutchings16} proves \eqref{eq:C} when $1\leq a\leq2$ and $\beta$ is a positive integer. Digiosia--Nelson--Ning--Weiler--Yang \cite{DNNWY22} establish additional cases for half-integer axis ratios and specified ranges of $a$. The new cases in Theorem \ref{thm:C} include non-integer ellipsoids in the target. Moreover, we emphasize that the role of \eqref{eq:polydisk-rigidity} is specifically to deal with the cases in (i) in Theorem \ref{thm:C} when $\beta \in (1,2)$. 

\begin{remark} There is a limitation to the inclusion criterion in Theorem \ref{thm:C}. For $\beta\geq2$, the region where $c>2$ and $\beta<a<2\beta+2$ is not covered by the argument above. As a matter of fact, for sufficiently thin polydisks, folding can improve on inclusion; compare \cite{Schlenk05} and \cite[Remark~1.12]{HZ26}. The standard one-fold construction $P(1,a) \hookrightarrow P(2, 1 + \frac{a}{2})$ yields
\[
 c_\beta(a)\leq2+\frac{a+2}{2\beta}.
\]
This is strictly smaller than $1+a/\beta$ precisely when $a>2\beta+2$. The line $a=2\beta+2$ is where these two upper bounds agree. \end{remark}

\begin{remark} Without applying the polydisk shape invariant (from Theorem \ref{cor:polydisk-rigidity} or Example \ref{ex:ellipsoid-polydisk}), one can still prove Theorem~\ref{thm:C} when the parameters lie in certain ranges (in particular, for non-integer ellipsoids $E(c, \beta c)$). For instance, Theorem~\ref{thm:D} below relies only on the Lagrangian reduced shape and can confirm the conclusion of Theorem~\ref{thm:C} when $2 \leq a \leq \beta$ (and for all $c \geq 1$). The difference is that the parameter range $1 \leq a \leq 2$ will be missing. This is where the polydisk shape invariant can provide more information. \end{remark} 


\begin{mainthm}\label{thm:D}
Let $\Omega\subset\R^2_{\geq0}$ be a bounded, relatively open, star-shaped domain in $\R_{\geq 0}^2$.  For any $\beta>0$, denote
\[
 M_\Omega(\beta)\coloneqq\max_{(u,v)\in\overline{\Omega}}
 \left(u+\frac v\beta\right).
\]
Then the following statements hold.
\begin{enumerate}[(i)]
\item Suppose that $\beta\geq2$ and that $M_\Omega(\beta)$ is attained at a point $(r_*,s_*)$ satisfying $r_*>0$ and $2r_*\leq s_*\leq\beta r_*$. Then for every $c>0$,
\[
 X_\Omega\into E(c,\beta c)
 \quad\Longleftrightarrow\quad c\geq M_\Omega(\beta).
\]
\item Suppose that $M_\Omega(1)$ is attained at a point $(r_*,s_*)$ with $r_*>0$ and $r_* \leq s_* \leq 2r_*$. Then for every $R>0$,
\[
 X_\Omega\into B^4(R)
 \quad\Longleftrightarrow\quad R\geq M_\Omega(1)=r_* + s_*.
\]
\end{enumerate}
\end{mainthm}

The proof of Theorem~\ref{thm:D} is given in Section~\ref{sec:proof-D}. Note that in both (i) and (ii) in Theorem \ref{thm:D}, the sharp cases $c =  M_\Omega(\beta)$ and $R = M_\Omega(1)$ are the inclusions of moment regions. Theorem \ref{thm:D} generalizes the quadrilateral embedding results of Miranda--Salom\~ao--Trejos \cite[Theorem~3.8 and Corollary~3.9]{MST26}, where they only treat integer ellipsoids, including the ball, for specified quadrilateral domains. Here, our line testing criterion applies to arbitrary bounded, relatively open, star-shaped moment regions satisfying the stated maximizing-point condition and allows every real ellipsoid ratio $\beta\geq2$, in particular non-integer ratios.

\begin{figure}[!htbp]
  \centering
  \includegraphics[scale=0.8]{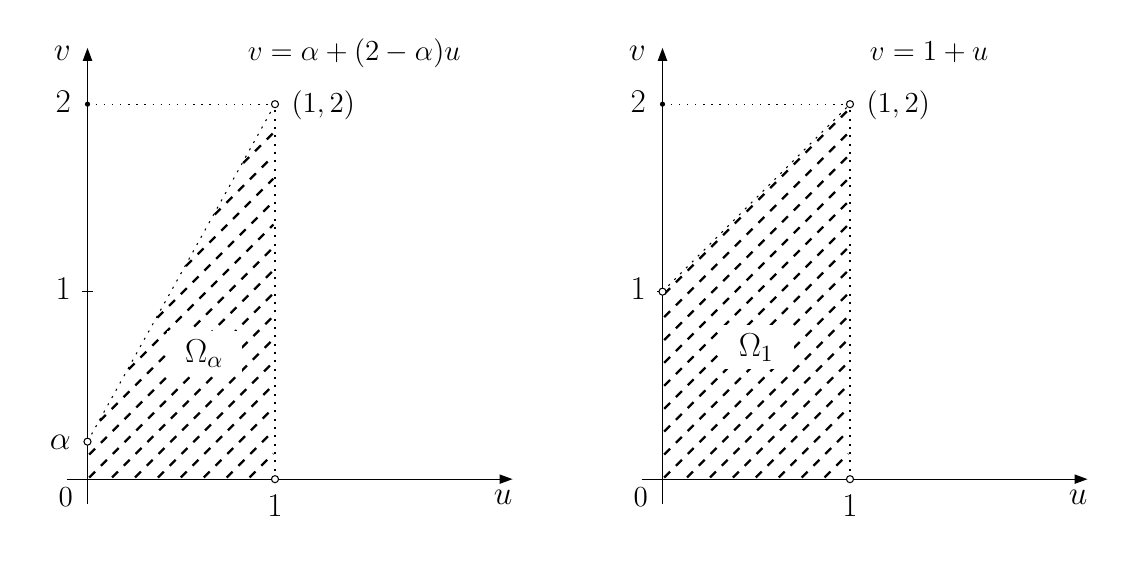}
  \caption{The moment regions of Example~\ref{ex:quadrilateral-family}
  and Corollary~\ref{cor:trapezoid}: $\Omega_\alpha$ with $0<\alpha<1$
  on the left, and $\Omega_1$ on the right.}
  \label{fig:quadrilateral-moment-regions}
\end{figure}

\begin{example}[A family of quadrilateral domains]\label{ex:quadrilateral-family}
For $\alpha>0$, let
\begin{equation} \label{mi-alpha}
 \Omega_\alpha\coloneqq\{(u,v)\in\R^2_{\geq0}\mid u<1,\ v<\alpha+(2-\alpha)u\},
 \qquad X_\alpha\coloneqq\mu^{-1}(\Omega_\alpha).
\end{equation}
The closure of $\Omega_\alpha$ has vertices $(0,0),(1,0),(1,2),(0,\alpha)$ (see the left picture in Figure \ref{fig:quadrilateral-moment-regions}). For every real $\beta\geq2$ and $0<\alpha\leq\beta+2$, the maximum of $u+v/\beta$ is $1+2/\beta$, attained at $(1,2)$. Then (i) of Theorem~\ref{thm:D} gives
\[
 X_\alpha\into E(c,\beta c)
 \quad\Longleftrightarrow\quad c\geq1+\frac2\beta.
\]
Likewise, for $0<\alpha\leq3$, part~\textup{(ii)} gives
\[
 X_\alpha\into B^4(R)\quad\Longleftrightarrow\quad R\geq3.
\]
Thus the criterion also enlarges the source-parameter range in \cite[Corollary~3.9]{MST26}, which assumes additionally $1\leq\alpha\leq2$. For $0<\alpha<2$, one can view $X_{\alpha}$ as a certain truncation of $P(1, 2)$ and, compared with Theorem~\ref{thm:C}, the embedding obstruction given in this example says that the obstruction won't be improved (though the volume drops) if we change $P(1,2)$ to $X_{\alpha}$. 
\end{example}

Furthermore, taking $\alpha=1$ (then see the right picture in Figure \ref{fig:quadrilateral-moment-regions} for $X_1$) and $\beta=3$ gives the following more concrete example.

\begin{cor}\label{cor:trapezoid} In terms of the notation in (\ref{mi-alpha}), for every $\lambda,R>0$,
\[
 \begin{aligned}
 X_{1}\into\lambda E(1,3) \quad\Longleftrightarrow\quad \lambda \geq \frac{5}{3}. 
 \end{aligned}
\]
However, $c_k^{\mathrm{ECH}}(X_1)=c_k^{\mathrm{ECH}}(E(1,3))$ for any $k \in \N$. \end{cor}

The detailed verification for the equalities of ECH capacities will be provided in Section~\ref{sec:proof-D}.

\begin{remark} For $X_1$ and $E(1,3)$, one should notice that their volumes and Gromov widths also coincide. Moreover, Ekeland--Hofer capacities help to distinguish them (see \cite{EH90, GH-capacity}), but give only $\lambda\geq4/3$. Therefore, the reduced shape supplies information missed by both capacity comparisons.\end{remark}

\subsection{Shape of intersections}\label{sec:intro-intersections}

By definition, the monotonicity of the Hamiltonian shape invariant gives an obvious inclusion
\begin{equation} \label{int-sh-inclusion}
 \Sh_H^+\!\left(\bigcap_{i=1}^n X_i\right)
 \subset\bigcap_{i=1}^n\Sh_H^+(X_i),
\end{equation}
and the same inclusion holds for $\Sh^+$. The following elementary example immediately shows that in general the equality in (\ref{int-sh-inclusion}) fails.

\begin{example}\label{ex:intersection-polydisks}
For $M>2$, let $X_1\coloneqq P(1,M)$ and $X_2\coloneqq P(2,2)$. Their intersection is the fixed polydisk $P(1,2)$. However, we have
\[
 \bigl(\Sh_H^+(X_1)\cap\Sh_H^+(X_2)\bigr)
 \setminus\Sh_H^+(X_1\cap X_2)
 =[1/2,1)\times[2,M).
\]
This rectangle lies in $\Wred$, so it also serves as the gap for the reduced shape invariant $\Sh^+$. Its area is $(M-2)/2$, which is unbounded as $M\to\infty$. \end{example}

The intersection theorem (more precisely, Hind--Kerman's computation on graphical type domains \cite[Theorem~1.17]{HK26}), our Theorems~\ref{thm:A} and \ref{thm:B} help to determine some cases where equality of (\ref{int-sh-inclusion}) holds. 

\begin{mainthm}\label{thm:E}
Let $n\geq1$ and let $f_i\colon[0,a_i]\to[a_i,\infty)$ be continuous, where $a_i>0$ for $1\leq i\leq n$. Then
\begin{equation}\label{eq:intersection-discrepancy}
 \begin{aligned}
 &\left(\bigcap_{i=1}^n \Sh_H^+(X_{\Omega_{f_i}})\right)\setminus\Sh_H^+\!\left(\bigcap_{i=1}^n X_{\Omega_{f_i}}\right) =\left\{(r,s)\;\middle|\;
 \begin{aligned}
 &\tfrac12\min_i a_i\leq r<\min_i a_i,\\
 &\min_i f_i(r)\leq s<\min_{\{i\mid a_i\leq2r\}}f_i(r)
 \end{aligned}\right\}.
 \end{aligned}
\end{equation}
Moreover,
\[
 \Sh^+(Y)=\Sh_H^+(Y)\cap\Cred
 \qquad\text{for }Y=X_{\Omega_{f_1}},\ldots, X_{\Omega_{f_n}},\ \bigcap_iX_{\Omega_{f_i}}.
\]
Therefore, if $\min_i f_i(r)= \min_{\{i\mid a_i\leq2r\}}f_i(r)$ for any $\tfrac12\min_i a_i\leq r<\min_i a_i$, then taking $\Sh_H^+$ commutes with the intersection; the same holds for $\Sh^+$ when restricted to $\Cred$.
\end{mainthm}

Theorem~\ref{thm:E} gives sharp embeddings from toric convex domains into toric concave domains. Recall that embeddings from toric concave domains to toric convex domains are completely determined by the ECH capacities. 

\begin{figure}[!htbp]
  \centering
  \includegraphics[scale=0.7]{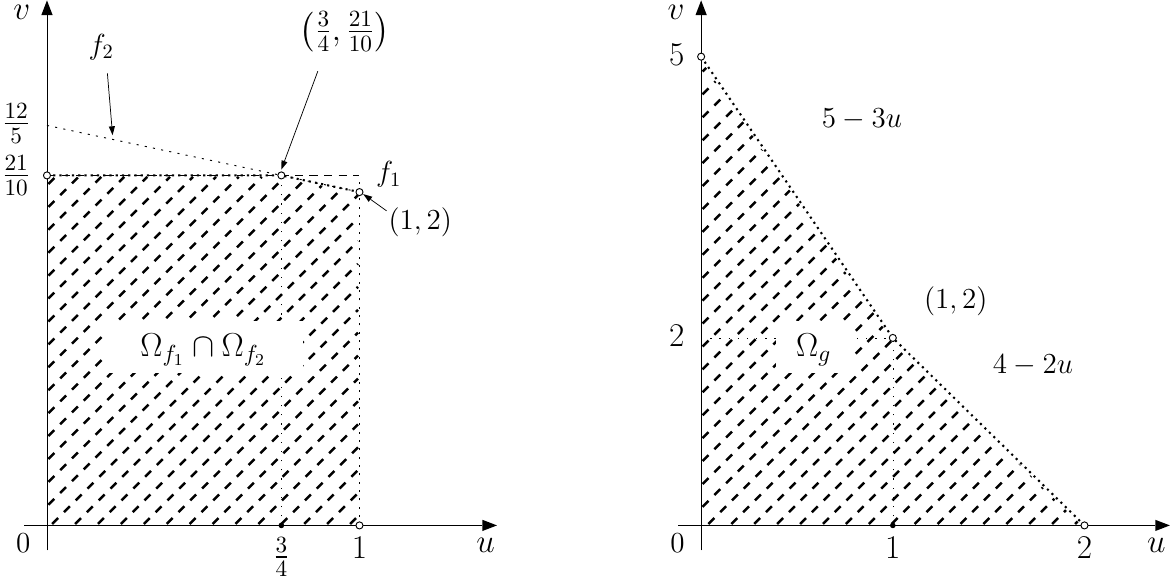}
  \caption{The moment regions in Corollary~\ref{cor:replacement-1-13}.}
  \label{fig:corollary-1-16-moment-regions}
\end{figure}

\begin{cor}\label{cor:replacement-1-13}
Consider the following two profile functions
\[
 f_1(u)=\frac{21}{10} \,\,\,\,\mbox{and} \,\,\,\,
 f_2(u)=\frac{12}{5}-\frac25u \,\,\,\,\,\,\,\,\mbox{$(0\leq u \leq 1)$}
\]
and denote $S:=X_{\Omega_{f_1}}\cap X_{\Omega_{f_2}}$; see the left picture in Figure \ref{fig:corollary-1-16-moment-regions}. Then we have the following equality $\Sh_H^+(S)=\Sh_H^+(X_{\Omega_{f_1}})\cap\Sh_H^+(X_{\Omega_{f_2}}).$ Moreover, let $T:=X_{\Omega_g}$ be the toric domain determined by the profile function 
\[
 g(u)=\max\{5-3u,4-2u\}\qquad (0\le u\le2);
\]
see the right picture in Figure \ref{fig:corollary-1-16-moment-regions}. Then for every $\lambda>0$, we have $S\into\lambda T$ if and only if $\lambda \geq 1$. 
\end{cor}

\begin{remark}\label{rmk-MST} For domains $S$ and $T$ in Corollary \ref{cor:replacement-1-13} above, a direct comparison of ECH capacities yields $\cECH{k}(S) \le\frac{31}{40}\cECH{k}(T)$ for all $k$. Moreover, Miranda--Salom\~ao--Trejos \cite{MST26} gives an embedding obstruction from toric convex domains into toric concave domains, which incorporates conditions beyond ECH capacities. Applied to the embedding $S \hookrightarrow \lambda T$, it can be verified, with the assistance of an LLM, that all the necessary conditions in \cite[Theorem~3.1]{MST26} are satisfied for $\lambda = \frac{19}{20}$. Nevertheless, Corollary \ref{cor:replacement-1-13} shows that $S \not\hookrightarrow \frac{19}{20}T$. Thus for this pair of domains, the shape invariant provides a strictly stronger obstruction than ECH capacities or the criterion \cite[Theorem~3.1]{MST26}. \end{remark}

Another consequence of Theorem~\ref{thm:E} is its application to the intersection of ellipsoids. 

\begin{cor}\label{prop:aligned-intersections}
Let $n\geq1$, let $0<a_i\leq b_i<\infty$ with arbitrary real axis ratios $b_i/a_i\geq2$. Then 
\[  \Sh^+\left(\bigcap_{i=1}^n E(a_i,b_i)\right) = \left(\left(\bigcap_{i=1}^n \mu(E(a_i, b_i))\right) \cup \left\{0< r < \frac{\min_{i=1}^n a_i}{2}\right\}\right) \cap \mathcal C. \]
The same formula holds for the Hamiltonian shape invariant $\Sh_H^+$, with $\mathcal C$ replaced by $\{(r,s)\mid0<r\leq s\}$, whenever $b_i/a_i \in \N_{\geq 2}$ for all $1 \leq i \leq n$. 
\end{cor}

Note that when viewing $E(a,b)$ as the graphical type toric domain $X_{\Omega_g}$ where $g:[0,a] \to [0, \infty)$ is linear, it does not satisfy the condition in Theorem~\ref{thm:E} directly (since the minimum of $g$ is $0$ instead of $a$). To initiate the application of Theorem~\ref{thm:E} and obtain the desired rigidities on (Hamiltonian) shape invariant, we will need to consider $E(a,b) \cup P(a,a)$. The proof of Corollary~\ref{prop:aligned-intersections} appears in Section~\ref{sec:proof-E}.

\subsection{Banach--Mazur geometry}

Let $\mathcal S$ be the space of bounded star-shaped domains in $(\R^4, d\lambda_{\rm std})$, without a volume normalization. We equip it with the \emph{Banach--Mazur pseudometric}\footnote{Here, $\log$ denotes the natural logarithm and the rescaling is still following (\ref{eq:scaling}).}
\begin{equation}\label{eq:BM}
 \dBM(X,Y)\coloneqq\inf\left\{\log\lambda \geq 0\, \bigg| \, 
 \frac{1}{\lambda} X\into Y\into\lambda X\right\}.
\end{equation}
The two arrows in \eqref{eq:BM} are independent symplectic embeddings; no condition is imposed on their composition. The computation or estimation of $d_{\rm BM}$ is rather difficult in general. For a subset $\mathcal S'\subset\mathcal S$, put
\[
 \dBM(X,\mathcal S')\coloneqq\inf_{Y\in\mathcal S'}\dBM(X,Y).
\]
The monotonicity of shape invariant immediately turns embeddings in both directions into scaled inclusions between subsets in $\R^2$. This suggests studying the geometry of $(\mathcal S,\dBM)$ through shape invariants.

Several recent results show that this geometry is considerably richer than that of convex domains. Dardennes--Gutt--Ramos--Zhang \cite{DGRZ2025} construct a sequence of monotone toric domains (in dimension-4) whose distance from the subset of convex domains tends to infinity. Since smooth strictly monotone toric domains are dynamically convex, this gives a quantitative distinction between dynamical convexity and symplectic equivalence to a convex domain. Soon afterward, Cristofaro-Gardiner--Hind \cite{CGH2023} show that $(\mathcal S,\dBM)$ admits quasi-isometric embeddings of $(\R^N,\|\cdot\|_\infty)$ for every $N$. This considerably generalizes \cite{DGRZ2025} since their resulting toric domains are all monotone. In the opposite direction, the convex counterexample to Viterbo's conjecture constructed by Haim-Kislev--Ostrover \cite{HKO2026}, together with the capacity comparison for monotone toric domains in \cite{GHR2022}, gives convex domains at positive distance from the set of monotone toric domains. Writing $\mathcal M$ for the subclass of monotone toric domains, the convex example $X_0$ in \cite{HKO2026} gives the quantitative bound
\begin{equation} \label{eq:previous-separation}
 \dBM(X_0,\mathcal M)\geq
 \frac14\log\left(\frac{3+\sqrt5}{5}\right)\approx 0.01153.
\end{equation}
More explicitly, $X_0$ is the Lagrangian product of pentagons in \cite{HKO2026}, and if there exists $Y \in \mathcal M$ and $\lambda \geq 1$ such that $X_0 \hookrightarrow \lambda Y$ and $Y \hookrightarrow \lambda X_0$, then we have 
\begin{align*}
c_{\rm EHZ}(X_0) \leq c_{\rm Z}(X_0) & \leq \lambda c_{\rm Z}(Y) \\
& = \lambda w_{\rm G}(Y) \leq \lambda^2 w_G(X_0) \leq \lambda^2 \sqrt{2 \Vol(X_0)}. 
\end{align*}
Therefore, $\lambda \geq \left(\frac{c_{\rm EHZ}(X_0)}{\sqrt{2 \Vol(X_0)}}\right)^{1/2} \geq \left(\frac{3+\sqrt5}{5}\right)^{1/4}$, where the last inequality shows up in the proof of Theorem 1.3 in \cite{HKO2026}. Here, $c_{\rm EHZ}$ denotes the Ekeland--Hofer--Zehnder capacity on convex domain $X_0$, $c_{\rm Z}$ denotes the cylindrical capacity, and $w_G$ denotes the Gromov width; over the subclass $\mathcal M$, we have $c_{\rm Z} = w_G$ by Theorem 1.7 in \cite{GHR2022}. Then by definition, we obtain the inequality (\ref{eq:previous-separation}).  

\medskip

These results motivate both quantitative separation questions for natural subclasses with respect to the geometry of $d_{\rm BM}$. Let $\Tplus$ and $\Tminus$ denote the subclasses of toric convex and toric concave domains, respectively. Here toric convex means that the reflection of the moment image across the coordinate axes is convex; toric concave means that its complement in the non-negative quadrant is convex. For a graphical moment region, these conditions correspond to a non-increasing concave or convex profile function, respectively. 

\begin{mainthm}\label{thm:F}
There exists a bounded open convex domain $X\subset(\R^4, d\lambda_{\rm std})$ such that
\[
 \dBM(X, \mathcal T_+)\geq \frac{1}{4} \log \frac{72}{53} \approx 0.07659.
\]
Here, $\mathcal T_+$ is the subclass of toric convex domains. Moreover, there exists a bounded open convex domain $Y\subset(\R^4, d\lambda_{\rm std})$ such that 
\[ 
 \dBM(Y, \Tpm)\geq \frac{1}{4} \log \frac{432}{233+54\sqrt{6}} \approx 0.04195.
 \]
Here, $\Tpm\coloneqq\Tplus\cup\Tminus$ is the union of subclasses of toric convex and toric concave domains.
\end{mainthm}

Note that the lower bound from Theorem~\ref{thm:F} is essentially larger than the one from (\ref{eq:previous-separation}); however, the subclass $\Tpm$ is more restrictive than $\mathcal M$. Compared with \cite{HKO2026}, the proof of Theorem~\ref{thm:F} is from a completely different perspective, and it is based on a volume-minimization process with prescribed Hamiltonian shape invariant. More explicitly, for a subset $S\subset\{0<r\leq s\}$ (which could be even disconnected), define
\begin{equation}\label{eq:shape-volume}
 \begin{aligned}
 \nu_\pm(S)&\coloneqq\inf\left\{\Vol(Y)\mid Y\in\Tpm,\ S\subset\Sh_H^+(Y)\right\};
 \end{aligned}
\end{equation}
similarly, one defines $\nu_+(S)$ over just the subclass $\Tplus$. These functionals measure the least symplectic volume needed to accommodate a prescribed subset of area classes. By a simple rescaling argument, we obtain the following estimate.

\begin{prop}\label{prop:BM-volume-bound}
Let $X\subset\R^4$ be a bounded open domain star-shaped with respect to the origin. Suppose that $S\subset\Sh_H^+(X)$ and $\nu_\pm(S)>0$. Then
\begin{equation}\label{eq:BM-volume-bound}
 \dBM(X,\Tpm)\geq\frac14 \max\left\{\log\!\left(\frac{\nu_\pm(S)}{\Vol(X)}\right), \,0\right\}. 
\end{equation}
The same conclusion holds for $d_{\rm BM}(X, \mathcal T_+)$ with respect to $\nu_+(S)$. 
\end{prop}

The short proof is given in Section~\ref{sec:proof-F}. The factor $1/4$ comes from the two embeddings and the quadratic scaling of volume, and an effective lower bound of $\dBM$ is provided precisely when $\nu_\pm(S)>\Vol(X)$. To apply Proposition~\ref{prop:BM-volume-bound}, one needs to control $\nu_\pm(S)$ for a preferred choice of subset $S$. Curiously, even prescribing a single point already leads to a non-trivial optimization problem, based on various rigidities of the Hamiltonian shape invariant.

\begin{figure}[!htbp]
  \centering
  \includegraphics[scale=0.6]{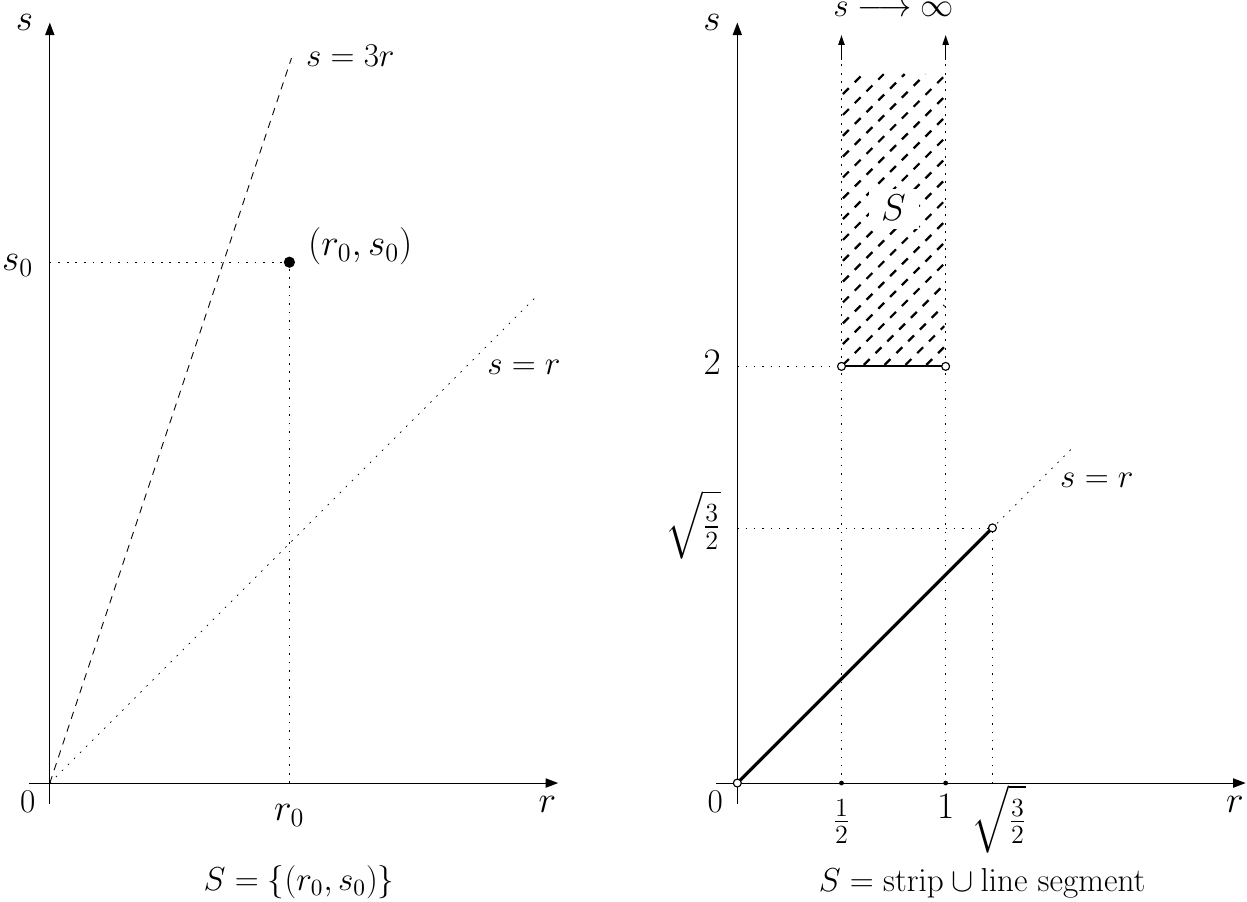}
  \caption{The prescribed subsets $S$ in Examples~\ref{ex:single-volume}
  and~\ref{ex:prescribed-shape}.}
  \label{fig:examples-1-19-1-20-prescribed-subsets}
\end{figure}

\begin{example}[One prescribed area pair]\label{ex:single-volume}
For any point $(r,s) \in \R^2_{>0}$ with $r \leq s$ (see the left picture in \ref{fig:examples-1-19-1-20-prescribed-subsets}), we have
\begin{equation}\label{eq:single-volume}
 \nu_+(\{(r,s)\})=\min\{rs,3r^2\}.
\end{equation}
The two terms have different geometric origins. The first is approached by polydisks containing the corresponding moment point. The second is related to the quadrilaterals $T_\rho$ from Theorem~\ref{thm:B}, whose areas are $3\rho^2$, as $\rho$ approaches $r$. A detailed proof is given in Section~\ref{sec:proof-F}.
\end{example}

For Theorem~\ref{thm:F}, we will consider a carefully chosen subset $S$, as demonstrated in the following example. Importantly, this subset $S$, though much more complicated than the one in Example~\ref{ex:single-volume}, still admits an exact computational value of $\nu_{\pm}$. 

\begin{example}[A chosen subset $S$ for Theorem~\ref{thm:F}]\label{ex:prescribed-shape}
Let
\begin{equation}\label{eq:S-configuration}
 S\coloneqq\left\{(r,s)\;\middle|\; \frac12<r<1,\ s\geq2\right\}
 \cup\left\{(r,r)\;\middle|\; 0<r<\sqrt{\frac32}\right\}.
\end{equation}
In other words, $S$ combines an unbounded region of non-monotone classes with a segment of diagonal classes (see the right picture in Figure \ref{fig:examples-1-19-1-20-prescribed-subsets}). Then we claim that 
\begin{equation} \label{s=3}
\nu_{\pm}(S) = 3. 
\end{equation}
Here, we give a complete argument. For any $\frac{2}{3}< r <1$, we know $(r, 3r) \in S$, therefore for any toric convex domain $Y \in \Tplus$ with $S \subset \Sh_H^+(Y)$, by Example~\ref{ex:single-volume} and definition (\ref{eq:shape-volume}) we have ${\rm Vol}(Y) \geq \min\{3r^2, 3r^2\} = 3r^2.$ Letting $r$ approach to $1$, we conclude that ${\rm Vol}(Y)  \geq 3$. On the other hand, for any toric concave domain $Y = X_{\Omega} \in \Tminus$, the result of Faisal--Li \cite[Corollary~12]{FL2026} identifies the Cieliebak-Mohnke capacity with the size of the diagonal of $\Omega$. In particular, due to the monotone component in our chosen $S$, for any $0<r<\sqrt{3/2}$ there exists a Lagrangian embedding from the product Lagrangian torus $L(r,r)$, and taking the supremum as $r$ approaches to $\sqrt{3/2}$ implies that the size of the diagonal of $\Omega$ must be at least $\sqrt{3/2}$. Then consider a line that passes through $(\sqrt{3/2}, \sqrt{3/2})$ such that the complement $\R_{\geq 0}^2 \backslash \Omega$ completely lies in one side of this line (where this line exists precisely due to the the geometry of a toric concave domain). Suppose the axis-intercepts are $A,B \geq 0$, then 
\[ {\rm Vol}(Y) \geq  AB/2 \geq 2(\sqrt{3/2})^2=3 \]
 since $\sqrt{3/2}/A + \sqrt{3/2}/B = 1$. This shows that $\nu_{\pm}(S) \geq 3$. Conversely, arbitrarily small enlargements of $Q_1$ from Theorem~\ref{thm:B} have their Hamiltonian shape invariants containing subset $S$, which implies $\nu_{\pm}(S) \leq 3$. This completes the proof of our claim (\ref{s=3}). 
\end{example}

\begin{remark} \label{rmk-on-nu-plus} The same argument in Example \ref{ex:prescribed-shape} shows that 
\[ \nu_+\left(\left\{(r,s)\;\middle|\; \frac12<r<1,\ s\geq2\right\}\right) = 3. \]
This will be used in the proof of the first conclusion in Theorem~\ref{thm:F}. The diagonal part $\{(r,r)\,|\, 0<r<\sqrt{3/2}\}$ in (\ref{eq:S-configuration}) is only for the consideration of the toric concave domain. Note that $\nu_-(\{(r,r)\,|\, 0<r<\sqrt{3/2}\}) = 3$ again by Faisal--Li \cite[Corollary~12]{FL2026}, and it can be realized by the symplectic ball $B^4(\sqrt{6})$.    \end{remark}

Towards the proof of Theorem~\ref{thm:F}, by the choice of subset $S$ in (\ref{eq:S-configuration}) as in Example \ref{ex:prescribed-shape}, we are left to construct a bounded open convex domain $X$ such that ${\rm Vol}(X)$ is sufficiently small (in particular, smaller than $3$). The construction is based on the proof of Theorem~\ref{thm:B}, tracing back to the wrapping argument in \cite[Theorem~5.2(2)]{HZ23}, and the resulting convex domain $X$ is a convex hull, where its volume can be calculated in a rather accurate way. 

\begin{remark} One can modify definition (\ref{eq:shape-volume}) by replacing ${\rm Vol}$ with some preferred symplectic capacity $c$. Accordingly, the lower bound of Proposition \ref{prop:BM-volume-bound} will be given by a ratio between $\nu_{\pm}(S)$ and $c(X)$. This could result in a better estimation, while the difficulty lies in the computation of $c(X)$, often not as straightforward as the volume.  \end{remark} 

\subsection{Local non-compactness} Due to the existence of unbounded strips in the (Hamiltonian) shape invariants, results on large-scale geometry (with respect to $d_{\rm BM}$) are difficult to obtain from this invariant, since under large rescaling those unbounded strips start to play the dominating roles and provide almost no obstruction at all. However, the following result shows that non-triviality of the local topology with respect to $d_{\rm BM}$ can be detected. Recall that $\mathcal S$ denotes the space of all bounded open star-shaped domains in $(\R^4, d\lambda_{\rm std})$.

\begin{mainthm}\label{thm:G}
The pseudo-metric space $(\mathcal S, d_{\rm BM})$ is {\rm not} locally compact. Explicitly, around the symplectic ellipsoid $E(1,3)$, for every $0<\ep<\log \frac{5}{3}$, there is a sequence of star-shaped toric domains $X_n\in\mathcal S$ with volumes converging to ${\rm Vol}(E(1,3)) = \frac{3}{2}$ such that
\begin{equation} \label{local-non-compactness}
 \dBM(X_n,E(1,3)) \leq \ep \,\,\,\,\mbox{and}\,\,\,\, \dBM(X_n,X_m) \geq \ep \,\,\,\,\mbox{for any $n\neq m$}. 
\end{equation}
Moreover, possibly using different bounds, we can choose our sequence of $X_n$ such that they all have the same volume. That is, the pseudo-metric $d_{\rm BM}$ remains locally non-compact when we restrict to the subset of $\mathcal S$ consisting of domains of a fixed volume.

The same conclusion holds for the quotient pseudo-metric $(\mathcal S/ \sim_{\mbox{\tiny{symp}}}, d_{\rm BM})$. \end{mainthm}

For the last conclusion in Theorem~\ref{thm:G}, the quotient space is defined by identifying two domains $X, Y \in \mathcal S$ if $Y = \Phi(X)$ for some symplectomorphism $\Phi$ in $(\R^4, d\lambda_{\rm std})$. Note that $d_{\rm BM}$ is well-defined in this quotient space. The construction of the desired sequence $X_n$ will be given in Section~\ref{sec:proof-G}. In particular, the embedding obstruction expressed in terms of (\ref{local-non-compactness}) is provided by Hind--Kerman's intersection theorem \cite[Corollary~1.8]{HK26}. 

\begin{remark} Unifying the volume of desired domains in the second conclusion in Theorem~\ref{thm:G} is to eliminate the possible influence of the volume on the estimation of $d_{\rm BM}$. Indeed, by rescaling we know that 
\[
 \dBM(U,V)\geq\frac12  \left|\log\frac{\Vol(U)}{\Vol(V)}\right|.
\]
Therefore, the obstruction in (\ref{local-non-compactness}) may simply come from the ratio difference of volumes, while the second conclusion in Theorem~\ref{thm:G} excludes this possibility.  
\end{remark}

\section*{Acknowledgement}
This work was completed during the second author's academic visit to the University of Notre Dame. He thanks the first author for his hospitality and Notre Dame for its active research environment. The first author is partially supported by the Simons Foundation. The second author is partially supported by the National Key R\&D Program of China (Grant No.~2023YFA1010500) and NSFC Grants No.~12301081, No.~12511540054, No.~12361141812, and No.~12631007.

\medskip

\noindent{{\bf About AI use:}} Besides the proofreading at the very last stage of writing up this paper and the assistance on figure drawing in a rather efficient way, LLMs were also involved in the following places: (i) verifying the necessary conditions in \cite[Theorem~3.1]{MST26} (see Remark \ref{rmk-MST}); (ii) in the proof of the second conclusion of Theorem~\ref{thm:F}, searching for the optimal deformation of the rectangle $R$ in (\ref{easy-R}) so that the resulting convex hull admits as small a volume as possible; (iii) computing the exact volume of the convex hull in (\ref{vol-conv-deform-R}) (see Remark \ref{rmk-LLM-compute-conv}).

\section{Proof of Theorem~\ref{thm:A}}\label{sec:proof-A}
\begin{proof} [Proof of Theorem~\ref{thm:A}]

We follow \cite{HK26} closely.



Let us fix some notation. Let $L=\varphi(\T^2)\subset Z^4(a)$, and let $e_1,e_2$ be the prescribed integral basis, with $\mu(e_1)=\mu(e_2)=2$ and $\omega(e_1)=r$, $\omega(e_2)=s$. First, note that every Maslov-two class has boundary $(1-k)e_1+ke_2$ and area
\[
 r+k(s-r),\qquad k\in\Z.
\]
For $k\le-1$ this is at most $2r-s\le0$; for $k=0$ it equals $r$; and for $k\ge1$ it is at least $s$. Consequently, a {\it positive} Maslov-two area is either $r$ or at least $s$, and the unique area-$r$ class has boundary $e_1$. Similarly, a Maslov-zero class has boundary $k(e_2-e_1)$, so its positive area is at least $s-r\ge r$.

\medskip

Suppose, to the contrary, that $L$ is disjoint from $L(r,t)$ for all $t\ge s$. Since $2r \ge a$ we can choose a $b$ sufficiently large that $L\subset P(2r,b)$. Compactify $P(2r, b)$ to $X:=S^2(2r)\times S^2(b)$, and denote the lines at infinity by $S_{\infty} := \{ \infty \} \times S^2$ and $T_{\infty} := S^2 \times \{ \infty \}$.

For a given compatible almost complex structure $J$, we have a foliation of $X$ by the images of $J$-holomorphic spheres representing the class $[S^2 \times \{*\}]$. The sphere $T_{\infty}$ is a leaf of the foliation, and each leaf intersects $S_{\infty}$ exactly once. Hence, as described in \cite[section 3.1.2]{HK26}, there is a projection map $p:X \to S_{\infty}$ mapping a point in $X$ to the intersection of the corresponding leaf with $S_{\infty}$.

We consider a family $J^N$ of such almost complex structures, stretched to length $N$ along $L \cup L(r,s)$. Taking a limit as $N \to \infty$ we obtain a finite energy $J^{\infty}$-holomorphic foliation of $X \setminus L$. 
Generic leaves of our finite energy foliation are (unbroken) holomorphic spheres, but there are also broken leaves, which are holomorphic buildings consisting of curves with asymptotic ends on $L \cup L(r,s)$ and of total area $2r$. In fact, exactly as in \cite[Lemma 3.4]{HK26}, these broken leaves consist of a pair of finite energy planes of Maslov class $2$, with ends matching a cylinder in either $T^* L$ or $T^* L(r,s)$. Exactly one of the planes avoids $S_{\infty}$, and by our area class assumption has area at least $r$, and if not $r$ exactly then at least $s$. But as $s \ge 2r$, the total area of our leaves, this plane must have area exactly $r$, and have boundary in the class $e_1$. The other plane in the broken leaf also has area $r$ and boundary now in the class $-e_1$.

The projection $p$ can be defined in the same way with respect to this limiting finite energy foliation, and extends continuously across $L$ and $L(r,s)$. The images of the two Lagrangians $p(L), p(L(r,s)) \subset S_{\infty}$ are disjoint embedded circles. Let $A$ be the component of $S_{\infty} \setminus p(L)$ which does not include $\{ \infty \}$, and let $C$ be the component of $S_{\infty} \setminus p(L(r,s))$ which does include $\{ \infty \}$.

The key continuity argument in \cite[section 3.2]{HK26}, see part (4) of Proposition 2.4 there, which based on the assumption that $L$ is disjoint form all of the $L(r,t)$, implies that $A$ and $C$ are disjoint. Then we can also define $B$ to be the cylinder $B:= S_{\infty} \setminus (\overline{A} \cup \overline{C})$.

Now we make an important observation: let $D$ be a Maslov $2$ disk (not necessarily holomorphic) with boundary on $L$ such that $p(D) = A$ and the projection $p:D \to A$ is injective, that is, $D$ meets each leaf of the foliation over $A$ exactly once. Suppose $\partial D$ lies in a class $ke_1 + le_2$. Then the Maslov $2$ condition implies
\[
k+l=1,
\]
and as $p|_D$ is injective $1$ we have
\[
l = \pm 1.
\]
Therefore $D$ has area either $s$ or $2r-s<0$.

Next, as in \cite[section 4]{HK26}, we fix $d+1$ points on $L(r,s)$ and $d$ points on $L$ and consider a limit of holomorphic spheres in the homology class $(1,d)$ under neck stretching. The result is a Type 3 building $\mathbf{F}$ consisting of the following:
\begin{itemize}
    \item an essential plane (that is, a curve on which $p$ has degree $1$) projecting to $A$ and having area $s$;
    \item an essential plane projecting to $C$ and having area $b-s$;
    \item an essential cylinder projecting to $B$ and having area $r$;
    \item $d$ broken planes asymptotic to $L(r,s)$;
    \item $d-1$ broken planes asymptotic to $L$.
\end{itemize}

Here we are using the observation above to conclude that the plane projecting to $A$ has area $s$ (it cannot have area $2r-s$ as $s \ge 2r$). As $L(r,s)$ is a product it is easy to check that the plane projecting to $C$ has area $b-s$, and we already know that the broken planes have area $r$. The cylinder projecting to $B$ then must have area $r$ simply because our building occupies a class $(1,d)$, which has area $2rb + 2b$, and the remaining curves have total area 
\[
s + (b-s) + dr + (d-1)r.
\]

By a similar argument we can fix $d$ points on $L(r,s)$ and $d+1$ points on $L$ and take limits of holomorphic spheres in the class $(1,d)$ to obtain a building $\mathbf{G}$ consisting of the following curves:
\begin{itemize}
    \item an essential plane (that is, a curve on which $p$ has degree $1$) projecting to $A$ and having area $s$;
    \item an essential plane projecting to $C$ and having area $b-s$;
    \item an essential cylinder projecting to $B$ and having area $r$;
    \item $d+1$ broken planes asymptotic to $L(r,s)$;
    \item $d$ broken planes asymptotic to $L$.
\end{itemize}

The buildings $\mathbf{F}$ and $\mathbf{G}$ can be perturbed to give two smooth symplectic spheres $F$ and $G$ which intersect positively and transversally at $d$ points in $p^{-1}(A)$ and at $d$ points in $p^{-1}(C)$. Let $H_1, \dots, H_{2d}$ be the fiber spheres through these intersection points.

Following \cite{HK26} again, we blow-up the $2d$ intersection points and then blow-down the proper transforms $\hat{H}_i$ of the $H_i$. The result is a copy $Y$ of $S^2 \times S^2$ in which the proper transforms $\hat{F}$ and $\hat{G}$ of $F$ and $G$ are disjoint and have self-intersection $0$. Let $T_0$ be a fiber sphere in $X$, disjoint from $F \cap G$, and with $p(T) \in A$. Then both of the proper transforms $\hat{T}_0$ and $\hat{T}_{\infty}$ are symplectic spheres which intersect each of $\hat{F}$ and $\hat{G}$ positively and transversally. It follows that
\[
Z := Y \setminus \left( \hat{F} \cup \hat{G} \cup \hat{T}_0 \cup \hat{T}_{\infty} \right)
\]
can be identified with a subset of $T^* {\T}^2$, actually of the form 
\[
\T^2 \times R \subset \T^2 \times \R^2 \equiv T^* \T^2
\]
where $R$ is a rectangle in $\R^2$.

By construction, the blowing up and down occurred in the complement of $L \cup L(r,s)$. Hence we can consider $L$ and $L(r,s)$ also as disjoint Lagrangian tori in $Z$. By \cite[Lemma 2.7]{HK26}, the tori are homologous to the zero section, and relatively exact.
Such an arrangement contradicts Gromov's Lagrangian intersection theorem \cite[Theorem in Section 2.3.$B_4''$]{Gromov85}.
\end{proof}

\section{Proof of Theorem~\ref{thm:B}}\label{sec:proof-B}

\begin{proof}[Proof of Theorem~\ref{thm:B}]
The quadrilateral $T_\rho$ is precisely $q(\rho,\rho)$ in
\cite[Definition~5.1]{HZ23}, so $Q_\rho=Q(\rho,\rho)$ in their notation.
Suppose first that $r<s$ and choose $\varepsilon>0$ with
$r+\varepsilon<\rho$. By \cite[Theorem~5.2(2)]{HZ23}, a Hamiltonian
image of $L(r,s)$ lies in
\[
 Q(r+\varepsilon,r+\varepsilon)
 =Q_{r+\varepsilon}\subset\operatorname{int}Q_\rho.
\]
If $s=r$, the product torus $L(r,r)$ itself lies in
$\operatorname{int}Q_\rho$. This proves the strip inclusion.
Moreover, if $X$ contains a neighborhood of the compact set $Q_\rho$,
then $Q_{\rho+\varepsilon}\subset X$ for sufficiently small
$\varepsilon>0$. Applying the same construction with $r=\rho<s$
gives the endpoint assertion; when $s=r=\rho$, the product torus
again suffices.
\end{proof}

\begin{proof}[Proof of Corollary~\ref{cor-strip-width}] Since $1<x<2$, we have a natural conclusion $E(1,x) \subset E(1,2)$. Meanwhile, we also have $E(1,x) \subset E(x,x)$. Therefore, 
\[ \Sh_H^+(E(1,x)) \subset \Sh_H^+(E(1,2)) \,\,\,\,\mbox{and}\,\,\,\, \Sh_H^+(E(1,x)) \subset \Sh_H^+(B^4(x)). \]
By Example \ref{ex:standard-shapes}, we know that $w(x) \leq \frac{1}{2}$ and $w(x) \leq \frac{x}{3}$. On the other hand, choose any $0<\rho<\frac{x}{x+2}$. Since $1<x<2$, the vertices of $T_\rho$ lie in $\mu(E(1,x))$:
\[ \rho+\frac{2\rho}{x}<1,\qquad 2\rho<1,\qquad 2\rho<x. \]
Thus $T_\rho\subset\mu(E(1,x))$, and Theorem~\ref{thm:B} gives $w(x)\geq\rho$. Taking the supremum over $\rho<\frac{x}{x+2}$ yields $w(x)\geq\frac{x}{x+2}$. This completes the proof. \end{proof}

We next derive the reduced shape invariant consequences of Theorems~\ref{thm:A} and~\ref{thm:B}.

\begin{proof}[Proof of Corollary~\ref{cor:graph-shape}]
Write $X\coloneqq X_{\Omega_f}$ and let $S$ denote the last set in
\eqref{eq:graph-shape}. Product tori realize
$\Omega_f \cap\Cred$ in $\Sh_H^+(X)$.
If $0<r<a/2$, choose $\rho$ with $r<\rho<a/2$. Then
$f(\rho)>2\rho$ and $f(2\rho)>0$; concavity and monotonicity imply $T_\rho\coloneqq\operatorname{conv}\{(0,0),(0,2\rho), (\rho,2\rho),(2\rho,0)\}  \subset\Omega_f.$
Therefore, Theorem~\ref{thm:B} realizes every pair $(r,s)\in\Cred$
with $r<a/2$. Consequently,
\[
 S\subset\Sh_H^+(X)\cap\Cred\subset\Sh^+(X).
\]

For the converse, let $L\subset X$ have reduced area pair $(r,s)$.
There is nothing to prove if $r<a/2$. Suppose first that
$r\geq a/2$ and $s\geq2r$. Since
$X\subset Z^4(a)\subset Z^4(2r)$, Theorem~\ref{thm:A}, applied
with cylinder parameter $2r$, gives
\[
 L\cap L(r,t_*)\neq\varnothing
 \qquad\text{for some }t_*\geq s.
\]
Then, under the moment map, the point $(r,t_*)\in\Omega_f$. Then since $\Omega_f$ is a graphical type moment image, for $s\leq t_*<f(r)$ we also have $(r,s)\in\Omega_f\subset S$.

It remains to treat the diagonal $s=r$ case. Set
\[
 \delta_f\coloneqq\sup\{u\in(0,a)\mid u<f(u)\}.
\]
The Lagrangian capacity computation for four-dimensional convex toric
domains \cite[Theorem~3.28]{Pereira22} gives $c_L(X)=\delta_f$;
the possibly non-smooth boundary can be handled by approximation.
Here $c_L$ is the Cieliebak--Mohnke capacity \cite{CM18}.
By definition, we know that the smallest positive symplectic area of a disk bounded by the monotone torus $L$ is $r$. Since $L$ is compact
and $X$ is open, $L\subset\tau X$ for some $0<\tau<1$.
Conformality gives
\[
 r\leq c_L(\tau X)=\tau\delta_f<\delta_f.
\]
As $f$ is non-increasing, this implies $(r,r)\in\Omega_f$. Hence, $\Sh^+(X)\subset S$, completing the proof.
\end{proof}

\begin{proof}[Proof of Corollary~\ref{cor:ellipsoid-shape}]
For $b>2a$, apply Corollary~\ref{cor:graph-shape} to the linear function $f(u)\coloneqq b(1-u/a)$, in particular, non-increasing and concave. In the restricted region $\Cred$, the reflected moment triangle contributes no additional points, since $r\leq s$ and $a\leq b$. Moreover, the case $b=2a$ is the integral ellipsoid, where its reduced shape invariant has been computed by \cite[Theorem~1.5]{HZ26}. This completes the proof. 
\end{proof}

\section{Proof of Theorem~\ref{thm:C}}\label{sec:proof-C}

\begin{proof}[Proof of Theorem~\ref{thm:C}]
The direction ``$\Longleftarrow$'' holds simply from the inclusion of
the moment images.  Indeed, if $a+\beta\leq\beta c$ and $(u,v)\in[0,1)\times[0,a)$,
then
\[
 u+\frac v\beta<1+\frac a\beta\leq c
\]
and $P(1,a) \subset E(c, \beta c)$. 

For the direction ``$\Longrightarrow$'' we work through the three cases.

For case~\textup{(i)}, suppose $1\leq c<2$ and there exists an embedding  $P(1,a)\into E(c,\beta c)$. Then after rescaling by $1/c$, this gives
\[
 \left(\frac1c,\frac{a}c\right)\in\Sh^P(E(1,\beta)) \,\,\,\,\mbox{where}\,\,\,\, \frac1 c>\frac12.
\]
Then by Example \ref{ex:ellipsoid-polydisk} (which comes from Theorem \ref{cor:polydisk-rigidity}), we know that point $\left(\frac1 c,\frac{a}c\right)$ lies in the moment image of the ellipsoid $E(1,\beta)$. In other words, we have
\[
 \frac1 c+\frac{a}{\beta c}<1.
\]
which is $a+\beta\leq\beta c$.

For case~\textup{(ii)}, the case $c<2$ is covered by the preceding argument.  So, let us assume that $c\geq2$ (and by assumption $a\leq\beta$), then
\[
 1+\frac a\beta\leq2\leq c \,\, (\mbox{which is equivalent to $a + \beta \leq \beta c$}). 
\]
This means for this range of parameters, the requested inequality holds automatically. 

For case~\textup{(iii)}, when $\beta =1$, the requested inequality, $1+a \leq c$, can be proved either by \cite[Theorem~1.5]{Hutchings16} or by a direct comparison between the Hamiltonian shape invariants of $B^4(c)$ and $P(1,a)$ (from \cite{HO20}). Moreover, when $c<2$, we get the inequality from case (i). The remaining case is when $\beta \geq 2$ and $c\geq 2$ (and by assumption $1\leq a \leq 2$), and then we have 
\[ a + \beta \leq 2 + \beta \leq 2\beta \leq \beta c. \]
Hence, for this range of parameters, the required inequality also holds. \end{proof}

\section{Proof of Theorem~\ref{thm:D}}\label{sec:proof-D}

\begin{proof}[Proof of Theorem~\ref{thm:D}]
For part~\textup{(i)}, suppose that $X_\Omega\into E(c,\beta c)$ and $c<M_\Omega(\beta)$. The hypotheses give
\[
 r_*<M_\Omega(\beta)=r_*+\frac{s_*}{\beta}\leq2r_*\leq s_*.
\]
Since $c<M_\Omega(\beta)$ and $2r_*\leq s_*\leq\beta r_*$, the radial segment from the origin to $(r_*,s_*)$ contains a point $(r,s)\in\overline\Omega$ satisfying
\[
 \frac c2<r<c<s, \qquad r+\frac{s}{\beta}>c.
\]
Since these inequalities are strict, they remain valid at a sufficiently nearby point still denoted by $(r,s)\in\Omega$. Then the monotonicity of the Hamiltonian shape invariant therefore gives a Hamiltonian image of $L(r,s)$ in $E(c,\beta c)\subset Z^4(c)$. By \cite[Corollary~1.8]{HK26}, after rescaling, this image must meet $L(r,w)$ for some $w\geq s$. However, every such torus lies outside the ellipsoid $E(c, \beta c)$, because
\[
 r+\frac w\beta\geq r+\frac s\beta>c.
\]
This is a contradiction.

For part~\textup{(ii)}, suppose that $X_\Omega\into B^4(R)$ with $R<r_*+s_*$. Since $r_*\leq s_*\leq2r_*$, we have
\[
 \min\{r_*,s_*\}=r_* \geq\frac{r_*+s_*}{3}>\frac R3.
\]
These inequalities are strict where needed, so a small perturbation gives a point $(r,s)\in\Omega$ satisfying
\[
 \min\{r,s\}>\frac R3,\qquad r+s>R.
\]
By the ball shape formula in Example~\ref{ex:standard-shapes}, no Hamiltonian image of $L(r,s)$ lies in $B^4(R)$. However, the assumed embedding and monotonicity of the Hamiltonian shape invariant imply that such an image exists, which provides the desired contradiction.\end{proof}

\begin{proof}[Proof of Corollary~\ref{cor:trapezoid}]
The moment region $\Omega_1$ satisfies the hypothesis of case~(i) in Theorem~\ref{thm:D}, where the function $u+v/3$ attains its maximum on the closure of $\Omega_1$ at $(1,2)$ and the corresponding value $M_{\Omega_1}(3) = 5/3$. Then the desired conclusion is from case~(i) in Theorem~\ref{thm:D}. In what follows, we will verify the coincidence of the ECH capacities. 

The closure of $X_{\Omega_1}$ is the quadrilateral domain
$Q_{1,1,1,2}$ in the notation of \cite{MST26}. Then their capacity formula \cite[Proposition~3.10]{MST26} gives
\begin{equation}\label{eq:trapezoid-ech-min}
 c_k^{\mathrm{ECH}}(X_{\Omega_1})
 =\min\left\{2m+n\;\middle|\;
 \begin{array}{l}
 m,n\in\Z_{\geq0},\\[2pt]
 (m+1)(n+1)+\dfrac{m(m+1)}2\geq k+1
 \end{array}\right\}.
\end{equation}
Although the cited proposition uses the closed domain, the formula remains the same for its interior. Next, we convert the minimum in \eqref{eq:trapezoid-ech-min} into a counting
problem. For $L\in\Z_{\geq0}$, by definition $c_k^{\mathrm{ECH}}(X_{\Omega_1})\leq L$ if and only if there exists some $(m, n) \in \Z_{\geq 0}^2$ satisfying $2m+n \leq L$ and $(m+1)(n+1)+{m(m+1)}/2 \geq k+1$. Introduce the following function of $L$, 
\begin{equation} \label{kappa-count}
\kappa(L)\coloneqq\max_{\substack{m,n\in\Z_{\geq0}\\2m+n\leq L}}
 \left((m+1)(n+1)+\frac{m(m+1)}2\right).
\end{equation}
Then we can view $c_k^{\mathrm{ECH}}(X_{\Omega_1})$ from a different perspective, in terms of $\kappa(L)$. Explicitly, 
\[
 c_k^{\mathrm{ECH}}(X_{\Omega_1}) \leq L \quad\Longleftrightarrow\quad \kappa(L) \geq k+1.
\]

Next, we investigate how to maximize the expression $(m+1)(n+1)+{m(m+1)}/2$ in (\ref{kappa-count}) in the specified region for $(m, n)$.  For a fixed $m$ with $0\leq2m\leq L$, this expression being maximized is
increasing in $n$, so its maximum occurs at $n=L-2m$. Then 
\[
 (m+1)(n+1)+\frac{m(m+1)}2= (m+1)\left(L+1-\frac32m\right).
\]
Consider the function $f(x) := (x+1) (L +1 - \frac{3}{2} x)$ (for $x\in \R$). It has two roots: $-1$ and $\frac{2}{3}(L+1)$. Therefore, the maximum is obtained at $\frac{L}{3} - \frac{1}{6}$. To reduce to the discrete case, write $L=3q+s$, where $q\in\Z_{\geq0}$ and $s\in\{0,1,2\}$. Then 
\[ \frac{L}{3} - \frac{1}{6} = q + \frac{s}{3} - \frac{1}{6}. \]
Since $|\frac{s}{3}-\frac{1}{6}|\leq\frac{1}{2}$, the integer $m=q$ is a maximizer (when $s=2$, $m=q+1$ is also a maximizer). Hence
\begin{align*}
\kappa(L) &= (q+1)  \left(L+1-\frac{3q}{2} \right) \\
& = (q+1) \left( 3q+s +1 - \frac{3q}{2} \right) = (q+1) \left(\frac{3q}{2}+s+1 \right).
\end{align*}

On the other hand, the  formula
\cite[Proposition~1.2]{Hutchings11} says that the capacities of $E(1,3)$
are the numbers $i+3j$, for $(i,j)\in\Z_{\geq0}^2$, arranged in
nondecreasing order with repetitions, starting at zero.
Therefore, the number of entries\footnote{An alternative way to argue is to work out an ECH capacities formula of $E(1,3)$, similarly to the one in (\ref{eq:trapezoid-ech-min}). Indeed, eventually by counting lattice points one gets 
\[ c_k^{\mathrm{ECH}}(E(1,3))
 =\min\left\{3m+n\;\middle|\;
 \begin{array}{l}
 m,n\in\Z_{\geq0},\\[2pt]
 (m+1)(n+1)+\dfrac{3m(m+1)}2\geq k+1
 \end{array}\right\}.\]
Then by the same argument via the $\kappa$-function of $L$ in (\ref{kappa-count}) but with the constraint that $3m+n \leq L$, we know that the maximum of this $\kappa$-function in this case is also equal to $(q+1) \left(\frac{3q}{2}+s+1 \right)$.}
  at most $L$ is
\begin{align*}
\#\{(i,j)\in\Z_{\geq0}^2\mid i+3j\leq L\} &=\sum_{j=0}^{q}(L-3j+1)\\
 &=(q+1)(L+1)-\frac32q(q+1) =(q+1)\left(\frac{3q}{2}+s+1\right).
\end{align*}
Consequently, for every $k,L\in\Z_{\geq0}$, we have shown that $c_k^{\mathrm{ECH}}(X_{\Omega_1})\leq L$ if and only if $c_k^{\mathrm{ECH}}(E(1,3))\leq L$. 
Since ECH capacities of both $X_{\Omega_1}$ and $E(1,3)$ take non-negative integer values, these
equivalences prove the desired coincidence at each $k$. This completes the proof. 
\end{proof}

\section{Proof of Theorem~\ref{thm:E}}\label{sec:proof-E}

\begin{proof}[Proof of Theorem~\ref{thm:E}]
By Hind--Kerman's graphical toric domain computation
\cite[Theorem~1.17]{HK26}, after rescaling, for every continuous function $g\colon[0,b]\to[b,\infty)$ we have 
\begin{equation}\label{eq:single-hamiltonian}
 \Sh_H^+(X_{\Omega_g})=\left(\Omega_g \cup \left\{0<r < \frac{b}{2} \right\} \right) \cap \{r \leq s\}.
\end{equation}

We first show that
\begin{equation}\label{eq:single-reduced}
 \Sh^+(X_{\Omega_g})
 =\Sh_H^+(X_{\Omega_g})\cap\Cred.
\end{equation}
The inclusion from the right hand-side to the left hand-side follows from the definition. Conversely, suppose $L\subset X_{\Omega_g}$ is an embedded Lagrangian torus with area pair $(r,s)$. The displacement energy implies that $r<b$. If $r<b/2$, then $(r,s)$ lies in the unbounded strip $\{0< r < b/2\}$. If $s=r$, then
$s=r<b\le g(r)$, so $(r,s)\in\Omega_g$. Therefore, we are left with the following case:
\[
 b/2\le r<b \,\,\,\,\mbox{and}\,\,\,\, s\ge2r.
\]
Applying Theorem A to $L\subset Z^4(b)$, we obtain
\[
 L\cap\bigcup_{t\ge s}L(r,t)\ne\varnothing.
\]
Thus $(r,t)\in\Omega_g$ for some $t\ge s$, and consequently $s\le t<g(r)$. So, $(r,s)\in\Omega_g$. Together with \eqref{eq:single-hamiltonian}, we obtain \eqref{eq:single-reduced}.

Now, for the family of graphical toric domains as in our hypothesis, set
\begin{equation} \label{min-function}
 a:=\min_{1\le i\le n}a_i \,\,\,\,\mbox{and}\,\,\,\, 
 f(u):=\min_{1\le i\le n}f_i(u)\quad(0\le u\le a)
\end{equation}
Note that here each $f_i$ is restricted to $[0,a]$ when forming the ``minimum'' function $f(u)$.
This function $f$ is continuous and satisfies $f\ge a$ over $[0,a]$, and $\bigcap_{i=1}^n X_{\Omega_{f_i}}=X_{\Omega_f}$. Therefore~\eqref{eq:single-hamiltonian} and~\eqref{eq:single-reduced}
apply both to the intersection and to every $X_i$. In particular, this proves
\[
 \Sh^+(Y)=\Sh_H^+(Y)\cap\Cred
 \qquad\text{for }Y=X_{\Omega_{f_1}},\,\ldots,\,X_{\Omega_{f_n}}, \, \bigcap_{i=1}^n X_{\Omega_{f_i}}.
\]

Now, we get back to the first conclusion on the difference between the Hamiltonian shape invariant of the intersection and the intersection of Hamiltonian shape invariants. For $0<r<a/2$, all these
Hamiltonian shape invariants contain every $(r,s)$ with $s\ge r$, hence, no difference. Meanwhile, choosing $j$ with $a_j = a = \min_{1\le i\le n}a_i$, then (\ref{eq:single-hamiltonian}) applies to $X_{\Omega_{f_i}}$, and in particular we know that $\bigcap_{i=1}^n \Sh_H^+(X_{\Omega_{f_i}}) \subset \Sh_H^+(X_{\Omega_{f_j}}) \subset \{0< r< a\}$. Thus the desired discrepancy can occur only when
\[
 a/2\le r<a.
\]
For such $r$, for simplicity define
\begin{equation} \label{F-function}
 F(r):=\min_{\{i |  a_i\le2r\}}f_i(r).
\end{equation}
Note that the index set for taking minimum is non-empty since $j$ with $a_j = a$ lies inside. 

Since $r<a\le a_i$, the computation \eqref{eq:single-hamiltonian} says that the $i$-th factor imposes no upper bound on $s$ when $a_i>2r$, while it imposes precisely $s<f_i(r)$ when $a_i\le2r$. Therefore, for $a/2\le r<a$,
\begin{align*}
 (r,s)\in\bigcap_{i=1}^n\Sh_H^+(X_{\Omega_{f_i}})
 &\iff r\le s<F(r),\\
 (r,s)\in\Sh_H^+\left(\bigcap_{i=1}^n X_{\Omega_{f_i}} \right)
 &\iff r\le s<f(r),
\end{align*}
where $f$ is the minimizing function defined in (\ref{min-function}). Since $f(r)\ge a>r$, subtraction yields
\begin{equation}\label{eq:discrepancy}
 \left(\bigcap_{i=1}^n\Sh_H^+(X_i)\right)
 \setminus\Sh_H^+\left(\bigcap_{i=1}^n X_{\Omega_{f_i}} \right)
 =\left\{(r,s) \,\bigg| \,\frac a2\le r<a,\quad f(r)\le s<F(r)\right\}.
\end{equation}
This is the desired conclusion for the difference. 

Finally, by~\eqref{eq:single-reduced},
the reduced discrepancy is the intersection of the right-hand
side of~\eqref{eq:discrepancy} with $\Cred$. In fact, every point of the discrepancy already satisfies $s>2r$.
Indeed, choose $j$ attaining $f(r)=f_j(r)$. Since
$f(r)\le s<F(r)$, this index $j$ does not belong to the index set $\{i \,| \, a_i\le2r\}$. Therefore,
\[
 s\ge f_j(r)\ge a_j>2r.
\]
Thus the discrepancies of the reduced shape invariant and the Hamiltonian shape invariant are identical:
\[
 \left(\bigcap_{i=1}^n\Sh^+(X_i)\right)\setminus\Sh^+\left(\bigcap_{i=1}^n X_{\Omega_{f_i}} \right)
 =\left\{(r,s) \, \bigg| \, \frac a2\le r<a,\quad f(r)\le s<F(r)\right\}.
\]
This completes the proof. \end{proof}

\begin{proof}[Proof of Corollary~\ref{cor:replacement-1-13}]
First, let us compute the Hamiltonian shape invariant of $S$. The profile curve of $S$ on $[0,1]$ is
\[
 f(u)=\min\left\{\frac{21}{10},\frac{12}{5}-\frac25u\right\}.
\]
Apply Theorem~\ref{thm:E} with $a_1=a_2=1$, using the profiles on their original interval $[0,1]$. On this interval,
\[
 f_1\geq1,\qquad f_2\geq2>1.
\]
For every $1/2\leq r<1$, both indices satisfy $a_i\leq2r$, so the function in (\ref{F-function}) is $F(r)=\min\{f_1(r),f_2(r)\}=f(r)$. Therefore, the discrepancy in Theorem~\ref{thm:E} is empty, proving that
\begin{equation}\label{eq:replacement-commutation}
\Sh_H^+(S)=\Sh_H^+(X_{\Omega_{f_1}})\cap\Sh_H^+(X_{\Omega_{f_2}}).
\end{equation}
This is the first conclusion. Next, on $[0,1]$ we have $f\ge2$ and $f(u)\le\frac{12}{5}-\frac25u\le5-3u=g(u)$, which implies that 
\begin{equation}\label{eq:replacement-inclusions}
 P(1,2)\subset S\subset T.
\end{equation}
Therefore, for any $\lambda\ge1$, the inclusion provides the desired embedding $S \hookrightarrow \lambda T$. Conversely, suppose now that $0<\lambda<1$ and $S\into\lambda T$.
Choose $\lambda<t<1$. Since $L(t,2t)\subset P(1,2)\subset S$,
the monotonicity of the Hamiltonian shape invariant gives $(t,2t)\in\Sh_H^+(S) \subset \Sh_H^+(\lambda T)$. Moreover, the profile of $\lambda T$ is
$g_{\lambda}(u) = \max\{5\lambda-3u,4\lambda-2u\}$ on $[0,2\lambda]$.
To apply the shape computations, we replace $g_{\lambda}$ by the larger profile $\widetilde{g_{\lambda}}$, where
\[
 \widetilde g_\lambda(u)
 =\max\{2\lambda,5\lambda-3u,4\lambda-2u\},
 \qquad 0\le u\le2\lambda.
\]
This profile is continuous and takes values in $[2\lambda,\infty)$. Then \cite[Theorem~1.17]{HK26} yields
\[
 \Sh_H^+(\lambda T)
 \subset
 (\Omega_{\widetilde g_\lambda}\cup \{0<r < \lambda\}) \cap \{r\le s\}.
\]
We claim that $(t, 2t)$ lies outside this region.
Indeed, if $t\ge2\lambda$, then $(t,2t)$ lies outside of $\Omega_{\widetilde{g}_{\lambda}} \subset \{ r < 2 \lambda\}$. But if $\lambda < t<2\lambda$, then 
$\widetilde g_\lambda(t)=2\lambda<2t$, so $(t,2t)$ is again outside the
graphical region. In this way, we obtain a contradiction, and then $S\into\lambda T$ implies $\lambda \geq 1$.  \end{proof}

\begin{proof}[Proof of Corollary \ref{prop:aligned-intersections}] Write $E_i:=E(a_i,b_i)$ and $X:=\bigcap_{i=1}^n E_i$. By Corollary \ref{cor:ellipsoid-shape}, we know 
\begin{equation}\label{eq:single-ellipsoid}
 \Sh^+(E(a_i, b_i))=\bigl(\mu(E(a_i, b_i))\cup \{0<r<a_i/2\}\bigr)\cap\Cred.
\end{equation}
We will combine this formula with Theorem \ref{thm:E}, applied to suitable enlargements of $E(a_i, b_i)$.

First, denote $a: = \min_{i=1}^n a_i$, then $a_i\ge a$ and $b_i\ge2a_i\ge2a$ imply for any $1 \le i \le n$, 
\[
 E(a,2a)\subset E(a_i, b_i) \,\,\,\,\mbox{which implies that $E(a,2a)\subset \bigcap_{i=1}^n E(a_i, b_i)$}.
\]
Then the Hamiltonian shape invariant of $E(a,2a)$ from \cite[Theorem~1.8]{HZ26} (or Example \ref{ex:standard-shapes}) contains the unbounded strip $\{0<r< a/2\}$. Therefore
\begin{equation}\label{eq:hamiltonian-lower}
\left(\left(\bigcap_{i=1}^n \mu(E(a_i, b_i))\right)  \cup \{0<r< a/2\} \right) \cap \{0<r \leq s\} \subset\Sh_H^+\left(\bigcap_{i=1}^n E(a_i, b_i)\right).
\end{equation}
In particular,
\begin{equation}\label{eq:lower-chain}
\left(\left(\bigcap_{i=1}^n \mu(E(a_i, b_i))\right)  \cup \{0<r< a/2\} \right) \cap\Cred\subset\Sh^+\left(\bigcap_{i=1}^n E(a_i, b_i)\right).
\end{equation}

For the other containing direction, define
\[
 g_i(u):=\max\left\{a_i,\,b_i\left(1-\frac{u}{a_i}\right)\right\} \,\,\,\,(0\le u\le a_i)
\]
and set
\[
 Y_i=E(a_i, b_i)\cup P(a_i,a_i) \,\,\,\,\mbox{and}\,\,\,\, Y=\bigcap_iY_i.
\]
For each profile $g_i$, as usual the corresponding graphical moment image is denoted by $\Omega_{g_i}:=\{(u,v)\in\R_{\ge0}^2:u<a_i,\ v<g_i(u)\}$. 
Then $Y_i=\mu^{-1}(\Omega_{g_i})$, and each $g_i$ is continuous
with $g_i\ge a_i$. Thus the family $\{Y_i\}_{i=1}^n$ satisfies the
hypotheses of Theorem E. For $a/2\le r<a$, as in the proof of Theorem~\ref{thm:E}, introduce two functions:
\[
 g(r)=\min_i g_i(r),\qquad
 G(r)=\min_{\{i:a_i\le2r\}}g_i(r).
\]
Then Theorem E implies that 
\begin{equation}\label{eq:enlarged-discrepancy}
 \left(\bigcap_i\Sh^+(Y_i)\right)\setminus\Sh^+(Y)
 =\left\{(r,s)\in\Cred \, \bigg| \, \frac a2\le r<a,\quad g(r)\le s<G(r)\right\}.
\end{equation}

Now, suppose for a contradiction that there exists some $(r,s)\in\Sh^+(\bigcap_{i=1}^n E(a_i, b_i))$ but not in the left hand-side of (\ref{eq:lower-chain}). Then $r\ge a/2$ and $(r,s)\notin \bigcap_{i=1}^n \mu(E(a_i, b_i))$. Meanwhile, by the monotonicity of the reduced shape invariant, we know $(r,s)\in\Sh^+(E_i)$ for every $i$. Hence, choosing an index $i$ with $a_i=a$, we get
$r<a$. More generally, ~\eqref{eq:single-ellipsoid} shows that
\begin{equation}\label{eq:active-bounds}
 a_i\le2r \quad\Longrightarrow\quad s<b_i\left(1-\frac{r}{a_i}\right) \le g_i(r).
\end{equation}
Consequently,
\begin{equation}\label{eq:upper-G}
 s<G(r).
\end{equation}
Moreover, since $(r,s)\notin \bigcap_{i=1}^n \mu(E(a_i, b_i))$, there is an index $j$ such that $s\ge b_j\left(1-\frac{r}{a_j}\right)$. By~\eqref{eq:active-bounds}, this index is not
counted in the definition of $G(r)$ (i.e., $a_j>2r$). Using $b_j\ge2a_j$, we obtain
\[
b_j\left(1-\frac r{a_j}\right)
       \ge2a_j-2r>a_j.
\]
It follows that $g_j(r)=b_j\left(1-\frac r{a_j}\right)$, and hence
\begin{equation}\label{eq:lower-g}
 g(r)\le g_j(r)=b_j\left(1-\frac{r}{a_j}\right)\le s.
\end{equation}
Equations~\eqref{eq:upper-G} and~\eqref{eq:lower-g} mean $(r,s)$ lies in the discrepancy region~\eqref{eq:enlarged-discrepancy}, so $(r,s)\notin\Sh^+(Y)$. This contradicts $X\subset Y$ and
$(r,s)\in\Sh^+(X)$.

Finally, suppose that every ratio $b_i/a_i$ is an integer at least two. Then we can use the full Hamiltonian shape formula for each integral ellipsoid from Example~\ref{ex:standard-shapes}, and the Hamiltonian version of Theorem~\ref{thm:E}. The same argument then applies with $\Sh_H^+$ in place of $\Sh^+$ and $\{(r,s)\mid0<r\leq s\}$ in place of $\Cred$. 
\end{proof}

\section{Proof of Theorem~\ref{thm:F}}\label{sec:proof-F}

\begin{proof}[Proof of Proposition~\ref{prop:BM-volume-bound}]
Let $Y\in\Tpm$, and let $\lambda\geq1$ satisfy $\lambda^{-1}X\into Y\into\lambda X$. Rescaling the first embedding gives $X\into\lambda Y$, so the monotonicity of the Hamiltonian shape invariant yields $S\subset\Sh_H^+(\lambda Y)$.
Since $\lambda Y\in\Tpm$, the second embedding and volume scaling imply
\[
 \nu_\pm(S)\leq\Vol(\lambda Y)
 =\lambda^2\Vol(Y)
 \leq\lambda^4\Vol(X).
\]
Dividing by $\Vol(X)$ and taking logarithms give the desired lower bound for $\log\lambda$, since $\log\lambda\geq0$. Taking the infimum over all such $\lambda$ and $Y$ proves \eqref{eq:BM-volume-bound}.
\end{proof}

\begin{proof}[Proof of Example~\ref{ex:single-volume}]
We first prove the two upper bounds in \eqref{eq:single-volume}. For every $\varepsilon>0$, the product torus $L(r,s)$ lies in the convex toric domain $P(r+\varepsilon,s+\varepsilon)$. So, 
\[
 \nu_+(\{(r,s)\})\leq (r+\varepsilon)(s+\varepsilon)=rs + \ep(r+s + \ep).
\]
Letting $\ep \to 0$, we get $ \nu_+(\{(r,s)\}) \leq rs$. Moreover, for every $\rho>r$, Theorem~\ref{thm:B} guarantees a Hamiltonian image of $L(r,s)$ in $\operatorname{int}Q_\rho$, which is a toric convex domain of volume $3\rho^2$. Letting $\rho \to r$ yields $\nu_+(\{(r,s)\})\leq3r^2$. This concludes the upper bound of \eqref{eq:single-volume}.

For the other direction, let $Y=X_\Omega\in\Tplus$ such that $(r,s)\in\Sh_H^+(Y)$. Suppose first that $s=r$. A Hamiltonian image of the monotone Lagrangian torus $L(r,r)$ has smallest positive disk area $r$, so $r\leq c_{\rm CM}(Y)$, the  Cieliebak-Mohnke capacity defined in \cite{CM18}. By \cite[Theorem~3.28]{Pereira22} we have
\[
 c_{\rm CM}(Y)=\delta_\Omega \,\,\,\, \mbox{where $\delta_\Omega\coloneqq\sup\{t>0\mid(t,t)\in\Omega\}$}.
 \]
For every $0<t<\delta_\Omega$, since $\Omega$ is convex, we know $[0,t] \times [0,t] \subset\Omega$. Since $r\leq c_{\rm CM}(Y) = \delta_\Omega$, letting $t\nearrow\delta_\Omega$ gives
\[
 \Vol(Y)\geq\delta_\Omega^2\geq r^2 =\min\{rs,3r^2\}.
\]
Next suppose that $r<s<2r$. Define the axis intercepts
\[
 a_\Omega:=\sup\{u\geq0\mid(u,0)\in\Omega\},\qquad
 b_\Omega:=\sup\{v\geq0\mid(0,v)\in\Omega\}.
\]
If $\min\{a_\Omega,b_\Omega\}\leq s$, then, after interchanging the coordinate planes if necessary, $Y\subset Z^4(s)$. Since $s/2<r<s$, the scaled form of \cite[Corollary~1.8]{HK26} implies that a Hamiltonian image of $L(r,s)$ in $Y$ intersects $L(r,t)$ for some $t\geq s$. Consequently $(r,t)\in\Omega$, and the coordinatewise downward closure of a toric convex moment region gives $[0,r]\times[0,s]\subset\Omega$. Hence $\Vol(Y)\geq rs$.
If instead $a_\Omega>s$ and $b_\Omega>s$, for all $(u,v) \in \Omega$, by convexity we have a convex hull
\[
 \operatorname{conv}\{(0,0),(s,0),(u,v),(0,s)\}\subset\Omega.
\]
The area of the convex hull is at least $s(u+v)/2$. Thus $\Vol(Y)<rs$ would imply $u+v<2r$ for every $(u,v)\in\Omega$, so $Y\subset B^4(2r)$. But by Example~\ref{ex:standard-shapes} the bounds $r<s<2r$ imply that $(r,s) \notin \Sh_H^+(B^4(2r))$, giving a contradiction. Therefore $\Vol(Y)\geq rs=\min\{rs,3r^2\}$ also in this range.

It remains to consider $s\geq2r$, with the same definitions of $a_\Omega$ and $b_\Omega$. We divide the argument into two cases. 

Case I: either $a_\Omega\leq2r$ or $b_{\Omega} \leq 2r$. Then up to symmetry (which preserves the Hamiltonian shape invariant) we can assume that $a_\Omega\leq2r$, and $Y\subset Z^4(2r)$. Let $L\subset Y$ be a Hamiltonian image of $L(r,s)$. Theorem~\ref{thm:A}, or \cite[Corollary~1.8]{HK26}, applies, and implies that 
\[
 L\cap L(r,t)\neq\varnothing
 \qquad\text{for some }t\geq s.
\]
Thus $(r,t)\in\Omega$. In particular, $[0,r]\times[0,s]\subset\Omega$, which implies that $\Vol(Y)\geq rs$.

Case II: both $a_\Omega>2r$ and $b_\Omega>2r$. Both $(2r,0)$ and $(0,2r)$ then belong to $\Omega$. For any nonzero $p=(u,v)\in\Omega$, by the convexity of $\Omega$ again, we have 
\[
\operatorname{conv}\{(0,0),(2r,0),p, (0, 2r)\} \subset \Omega
\]
which implies that 
\[ 
\Vol(Y)\geq\operatorname{Area}(\operatorname{conv}\{(0,0),(2r,0),p, (0, 2r)\}) \geq r(u+v).
\]
If $\Vol(Y)<3r^2$, then $u+v<3r$ for every $(u,v)\in\Omega$. This means $Y\subset B^4(3r)$ and $\Sh_H^+(Y) \subset \Sh_H^+(B^4(3r))$. Our $(r,s) \in \Sh_H^+(Y)$ will also lie in $\Sh_H^+(B^4(3r))$. By the ball case in Example~\ref{ex:standard-shapes}, this is impossible, because otherwise $r<(3r)/3$ or $r+s<3r$ (which implies $s< 2r$). This contradiction then shows that $\Vol(Y)\geq3r^2 \geq \min\{rs,3r^2\}$ in this case.

In conclusion, we have shown that any $Y$ with $(r,s)\in\Sh_H^+(Y)$ has volume at least $\min\{rs,3r^2\}$. Taking the infimum over $\Tplus$ and combining with the upper bounds proves \eqref{eq:single-volume}.
\end{proof}

Now, let us proceed to the proof of Theorem~\ref{thm:F}. The proof essentially uses the wrapping argument from \cite[Theorem 5.2]{HZ23}, so we will illustrate how this works via the following lemma. For the proof, it will be convenient to work in Euclidean coordinates (instead of the complex or polar coordinates), so let us fix the notations first. Denote coordinates by $z_1=p+\sqrt{-1} q$ and $z_2=x+\sqrt{-1}y$, so 
\begin{equation} \label{omega-std}
 \omega_{\rm std}=dp\wedge dq+dx\wedge dy.
\end{equation}
Let $L(r,s)=\{\pi(p^2+q^2)=r,\ \pi(x^2+y^2)=s\}$. Consider the following polytope 
\begin{equation} \label{dfn-K} 
 K: =\left\{(p,q,x,y) \in \R_{\ge 0}^4 \,\bigg| \, \begin{array}{l} p\le1,\quad y\le1, \quad  q+y\le2\\
x\le1+p,\quad x\le1+q+y \end{array}\right\}.
\end{equation}
This polytope has $17$ vertices and $9$ three-dimensional facets. 

\begin{lemma} \label{lem:wrapping}
For every $b>1$ and every $\ep>0$, there exists a Hamiltonian diffeomorphism $\phi \in\Ham_c(\R^4,\omega_{\rm std})$ such that
\begin{equation}  \label{eq:wrapping}
 \phi(L(1,b))\subset \{(p, q, x, y) \in\R^4 \,| \, {\rm dist}_{\R^4}((p,q, x, y), K)<\ep\}.
\end{equation}
\end{lemma}

\begin{proof} Fix a sufficiently small $\eta>0$ such that $1+4\eta<b$. Choose a smooth planar disk
\[ A_\eta \subset (0,1)_p\times(0,1+\eta)_q \,\,\,\mbox{such that ${\rm area} (A_\eta)=1$}. \]
It exists because the containing rectangle has area \(1+\eta>1\). Let \(\Gamma_b\) be the boundary of a smooth rounded rectangle of area \(b\), arbitrarily close to the boundary of \([0,1]_x\times[0,b]_y\). Then the product \(\partial A_\eta\times\Gamma_b\) is a Hamiltonian image of the product torus \(L(1,b)\). Moreover, this Hamiltonian diffeomorphism can be realized by a compactly supported Hamiltonian isotopy after an appropriate cut-off (near its Hamiltonian trace).

Take a smooth nonincreasing function $\rho: \R \to [0,1+2\eta]$ where 
\[
 \rho(y)=1+2\eta\ \text{for $y \leq 0$},\qquad \rho(y)=0\ \text{for $y \geq 1+4\eta$},
\]
and with $\rho'\in(-1,0]$ for all $y$.
Such a choice is possible because \(1+2\eta<1+4\eta\). Now, consider the following Hamiltonian function, 
\[
 H_\eta(p,q,x,y)=\rho(y) \cdot p.
\]
With the convention \(\iota_{X_{H_{\eta}}}\omega_{\rm std}=-dH_{\eta}\), we obtain the Hamiltonian vector field as 
\[
 X_{H_\eta}=\rho(y)\partial_q-\rho'(y)p\partial_x.
\]
So, the resulting Hamiltonian diffeomorphism $\phi_{H_{\eta}}^1$ acts as follows, 
\begin{equation}  \label{eq:suspension}
 (p,q,x,y)\longmapsto
 \bigl(p,q+\rho(y),x-\rho'(y)p,y\bigr) =: \phi_{H_{\eta}}^1(p, q, x, y). 
\end{equation}

Since $-\rho'(y)\in[0,1)$, the displacement of the $x$-coordinate is nonnegative and at most $p$.
Thinking of our torus \(\partial A_\eta\times\Gamma_b\) as fibered over $\Gamma_b$, the disk $A_{\eta}$ over the bottom $\{y=0\}$ of the rectangle $[0,1]_x\times[0,b]_y$ is translated by \((0,1+2\eta)\) 
and so are displaced by $\phi_{H_{\eta}}^1$. 

Next, we will apply the wrapping map in the proof of Theorem 5.2 in \cite{HZ23}, it  keeps the already deformed lower part of the rectangle $[0,1]_x\times[0,b]_y$ fixed and folds the upper part of $[0,1]_x\times[0,b]_y$ into an arbitrarily small neighborhood of
\begin{equation} \label{deform-wrap}
 \overline A_\eta\times\left([0,1]_x\times[0,1+4\eta]_y\right).
\end{equation}
Importantly, since the wrapping, when viewed from the $xy$-plane, will pass through the bottom $[0,1]_x \times \{y=0\}$, in order to avoid self-intersection in $\R^4$, it was necessary to displace the disk $A_{\eta}$ from itself. In order to keep track of the trace of the deformation, we introduce some notation:
\begin{itemize}
\item{} $C_1:=[0,1]_p \times [0,1]_q \times [0,1]_x \times [0,1]_y$;
\item{} $T_0 :=\{(p,q+1-y, \lambda p,y) \,| \, 0\le p,q,\lambda,y\le1\}$;
\item{} $T_1 :=\{(p,q+1-y, 1+\lambda p,y) \,| \, 0\le p,q,\lambda,y\le1\}$;
\item{} $B: =\{(p,q+1,x,0) \,| \, 0\le p,q,x\le1\}$.
\end{itemize}
By (\ref{eq:suspension}), the deformed two vertical sides of the lower part of the rectangle $[0,1]_x\times[0,b]_y$, i.e., $\{x=0\} \times [0, 1+4\eta]$ and $\{x=1\} \times [0, 1+4\eta]$, lie within a neighborhood of \(T_0\cup T_1\); the bottom lies within a neighborhood of \(B\). Moreover, by (\ref{deform-wrap}) after wrapping the resulting domain (which takes care of the upper part of the rectangle $[0,1]_x\times[0,b]_y$, i.e., $[0,1]_x\times(1+4\eta,b]_y$) lies in a neighborhood of $C_1$. Observe that $B \subset {\rm Conv}(T_0 \cup T_1)$, the convex hull of $T_0$ and $T_1$. For the prescribed $\ep>0$, choose $\eta$ and the wrapping-neighborhood thickness sufficiently small. Then there exists a Hamiltonian diffeomorphism $\phi$ that maps $L(1,b)$ into the $\ep$-neighborhood of ${\rm Conv}(C_1 \cup T_0 \cup T_1 \cup B) = {\rm Conv}(C_1 \cup T_0 \cup T_1)$. 

Finally, we show that $K = {\rm Conv}(C_1 \cup T_0 \cup T_1)$. By construction, it is readily verified that $C_1, T_0, T_1 \subset K$, which implies that ${\rm Conv}(C_1 \cup T_0 \cup T_1) \subset K$ (since $K$ is convex). On the other hand, fix \(z=(p_\ast,q_\ast,x_\ast,y_\ast)\in K\), and put
\[
 t: =q_\ast+y_\ast\in[0,2].
\]
Then we will give an explicit convex decomposition in three cases.

\medskip

\noindent {\rm Case I}: \(1\le t\le2\).
Set
\[
 z_0=(p_\ast,q_\ast,0,y_\ast),\qquad
 z_1=(p_\ast,q_\ast,1+p_\ast,y_\ast).
\]
Putting \(p=p_\ast\), \(q=t-1\in[0,1]\), \(y=y_\ast\) and \(\lambda=0\) shows \(z_0\in T_0\), since then \(q+1-y_\ast=t-y_\ast=q_\ast\). Similarly putting \(p=p_\ast\), \(q=t-1 \), \(y=y_\ast\) and \(\lambda=1\) shows \(z_1\in T_1\).

Moreover, since \(0\le x_\ast\le1+p_\ast\) and \(1+p_\ast>0\), the number
\(\theta=x_\ast/(1+p_\ast)\) lies in \([0,1]\). Therefore
\begin{equation}
 z=(1-\theta)z_0+\theta z_1\in {\rm Conv}(C_1 \cup T_0 \cup T_1).
 \label{eq:case1}
\end{equation}

\smallskip
\noindent {\rm Case II}: \(0\le t<1\) and \(0\le x_\ast\le1\).
Since \(q_\ast,y_\ast\ge0\), one has \(q_\ast,y_\ast\le t<1\).
Together with \(0\le p_\ast,x_\ast\le1\), this gives
\(z\in C_1\subset {\rm Conv}(C_1 \cup T_0 \cup T_1)\).

\smallskip
\noindent {\rm Case III}: \(0\le t<1\) and \(x_\ast>1\).
Put $u:=x_\ast-1$. The inequalities \(x_\ast\le1+p_\ast\) and \(x_\ast\le1+t\) give $0<u\le p_\ast\le1$ and $u\le t<1$. Then \(t>0\) and \(1-u>0\) and we can define
\begin{align}
 a:=\left(1,\frac{q_\ast}{t},2,\frac{y_\ast}{t}\right) \,\,\,\,\mbox{and}\,\,\,\, c:=\left(
 \frac{p_\ast-u}{1-u},\
 \frac{q_\ast(t-u)}{t(1-u)},\
 1,\
 \frac{y_\ast(t-u)}{t(1-u)}
 \right).
\end{align}

First, \(a\in T_1\). Indeed, in the definition of $T_1$, we can choose
the parameters \(p=1\), \(q=0\), \(\lambda=1\), and
\(y=y_\ast/t\in[0,1]\).
Next, \(c\in C_1\). Indeed, its first coordinate satisfies $0\le\frac{p_\ast-u}{1-u}\le1$ by our choice of $u$. Its second and fourth coordinates are nonnegative,
and their sum is
\[
 \frac{q_\ast(t-u)+y_\ast(t-u)}{t(1-u)}
 =\frac{t-u}{1-u}\in[0,1],
\]
so these coordinates also lie in \([0,1]\). Finally, direct substitution gives
\begin{align*}
 u a+(1-u)c
 &=\left(
 u+p_\ast-u,\
 \frac{uq_\ast+q_\ast(t-u)}{t},\
 2u+1-u,\
 \frac{uy_\ast+y_\ast(t-u)}{t}
 \right)\\
 &=(p_\ast,q_\ast,1+u,y_\ast)
 =(p_\ast,q_\ast,x_\ast,y_\ast)=z.
\end{align*}
As \( u \in (0,1) \), we have
\(z\in {\rm Conv}(C_1\cup T_1)\subset {\rm Conv}(C_1 \cup T_0 \cup T_1)\),
completing the proof. 
\end{proof}

Now, we are ready to give the proof of the first conclusion of Theorem~\ref{thm:F}. 

\begin{figure}[!htbp]
  \centering
  \includegraphics[scale=0.65]{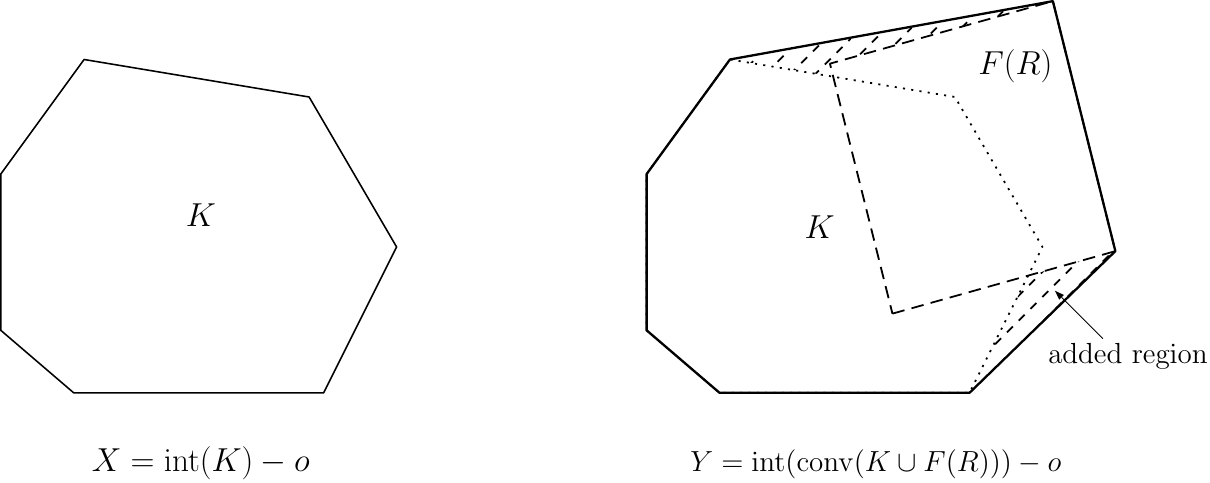}
  \caption{A schematic two-dimensional  picture of the construction in
  Theorem~\ref{thm:F}.}
  \label{fig:theorem-F-convex-domains}
\end{figure}

\begin{proof} [Proof of the first conclusion of Theorem~\ref{thm:F}] 

Suppose $(r,s) \in \{\frac{1}{2}<r<1, \, s\geq 2\}$. Then Lemma \ref{lem:wrapping}, together with a small rescaling argument, implies that $(r,s) \in \Sh_H^+({\rm int}(K))$ for the polytope $K$ constructed in (\ref{dfn-K}).

Next, set $o:=(1/2,1/2,1/2,1/2)\in{\rm int}(K)$ and $X:={\rm int}(K)-o$ (the left picture of Figure \ref{fig:theorem-F-convex-domains} shows a schematic 2-dimensional picture of this polytope). Then $X$ is bounded, open, convex, and contains the origin, and has the same shape invariants and volume as $K$. 
With $S$ defined by ~\eqref{eq:S-configuration}, by the second conclusion in Proposition~\ref{prop:BM-volume-bound} and Remark~\ref{rmk-on-nu-plus}, we obtain
\begin{equation} \label{est-nu-plus} 
d_{\rm BM}(X, \mathcal T_+) \geq\frac14 \max\left\{\log\!\left(\frac{\nu_+(S)}{\Vol(K)}\right), \,0\right\} = \frac14 \max\left\{\log\!\left(\frac{3}{\Vol(K)}\right), \,0\right\}. 
\end{equation}
What's left is the calculation of $\Vol(K)$. Dropping the inequality \(x\le 1+q+y\) from the definition of \(K\) gives the larger polytope
\[
V:=\left\{(p,q,x,y)\in\mathbb R^4 \,| \, 
0\le p,y\le1, \,
0\le q\le2-y,\,
0\le x\le1+p
\right\}.
\]
For fixed \(p,y\), the \((q,x)\)-section is a rectangle of area
\((2-y)(1+p)\). Therefore
\[
\operatorname{Vol}(V)
=\int_0^1\int_0^1(1+p)(2-y)\,dy\,dp
=\left(\int_0^1(1+p)\,dp\right)
 \left(\int_0^1(2-y)\,dy\right)
=\frac94.
\]
Meanwhile, the part removed by imposing \(x\le1+q+y\) satisfies $1+q+y<x\le1+p$, which is possible precisely when \(q+y<p\). Thus, up to boundary sets of volume zero, its coordinate ranges are
\[
0\le p\le1,\qquad
0\le y\le p,\qquad
0\le q\le p-y,\qquad
1+q+y<x\le1+p.
\]
The original bound \(q\le2-y\) is automatic here, since \(q\le p-y\le1-y\le2-y\). So, the length of the removed \(x\)-interval is \(p-q-y\), and 
\begin{align*}
\operatorname{Vol}(V\setminus K)
&=\int_0^1\int_0^p\int_0^{p-y}
  (p-q-y)\,dq\,dy\,dp\\
&=\int_0^1\int_0^p\frac{(p-y)^2}{2}\,dy\,dp =\int_0^1\frac{p^3}{6}\,dp =\frac1{24}.
\end{align*}
As a result, $\Vol(K) = \Vol(V) - \Vol(V \backslash K) = \frac94-\frac1{24} = \frac{53}{24} (<3)$. Therefore, the estimation (\ref{est-nu-plus}) yields the desired conclusion. 
\end{proof}

\FloatBarrier

The proof of the second conclusion of Theorem~\ref{thm:F} is based on the construction of $K$ in (\ref{dfn-K}). The desired convex domain is some {\it convex} enlargement of $K$ so that its Hamiltonian shape invariant also contains the diagonal part $\{(r,r) \,| \, 0<r<\sqrt{3/2}\}$ in (\ref{eq:S-configuration}). 
Thus $(r,r)\in\Sh_H^+({\rm int}(R))$.


\begin{proof}[Proof of the second conclusion of Theorem~\ref{thm:F}] Consider the following rectangle 
\begin{equation} \label{easy-R}
R: = [0,1]_p \times [0, \sqrt{3/2}]_q \times [0, \sqrt{3/2}]_x \times [0,1]_y  \subset \R^4. 
\end{equation}
For every $0<r<\sqrt{3/2}$, the interior of $R$ contains a Hamiltonian image of $L(r,r)$. Consider the following linear transform, as a shear map, 
\[ (p, q, x, y) \to \left(p, q - \frac{y}{2}, x + \frac{p}{2}, y \right) \]
where adding $-y/2$ in the second coordinate shifts $R$ in (\ref{easy-R}) toward smaller $q$, while adding $p/2$ in the third coordinate shifts it toward larger $x$. These directions are compatible with the sloping faces $q+y\le2$ and $x\le1+p$ of $K$ in (\ref{dfn-K}). Then add another shift by vector $(0, \frac{1}{3}, \frac{7}{6} - \sqrt{\frac{3}{2}}, 0)$. Together, we obtain an affine map $F: \R^4 \to \R^4$ defined by 
\begin{equation} \label{linear-F} 
F(p,q,x,y): =  \left(p,\ q-\frac y2+\frac13,\ x+\frac p2+\frac76-\sqrt{\frac{3}{2}},\ y\right).
\end{equation}
The corresponding convex hull ${\rm Conv}(K \cup F(R))$ has
\begin{equation} \label{vol-conv-deform-R}
\Vol({\rm Conv}(K \cup F(R))) = \frac{233+54\sqrt{6}}{144} \approx 2.53661 (<3). 
\end{equation}
Then set $o:=(1/2,1/2,1/2,1/2)$ and $Y:={\rm int}({\rm Conv}(K\cup F(R)))-o$ (the right picture of Figure \ref{fig:theorem-F-convex-domains} shows a schematic 2-dimensional picture of this convex hull). This is a bounded open convex domain containing the origin. 
By the first conclusion of Proposition \ref{prop:BM-volume-bound} and Example \ref{ex:prescribed-shape}, we know that 
\begin{align*} 
d_{\rm BM}(Y, \mathcal T_\pm) & \geq\frac14 \max\left\{\log\!\left(\frac{\nu_\pm(S)}{\Vol(Y)}\right), \,0\right\} \\
& = \frac14 \log\!\left(\frac{3}{\Vol(Y)}\right) = \frac14\log\!\left(\frac{432}{233+54\sqrt{6}}\right) 
\end{align*}
which is the required bound and completes the proof. 
\end{proof}

\begin{remark} \label{rmk-LLM-compute-conv} The calculation of the volume of convex hull (\ref{vol-conv-deform-R}) was done by an {\rm LLMs}. The calculation is similar to the one above, but a little more complicated.
\end{remark}

\section{Proof of Theorem~\ref{thm:G}}\label{sec:proof-G}

\begin{figure}[!htbp]
  \centering
  \includegraphics[scale=0.65]{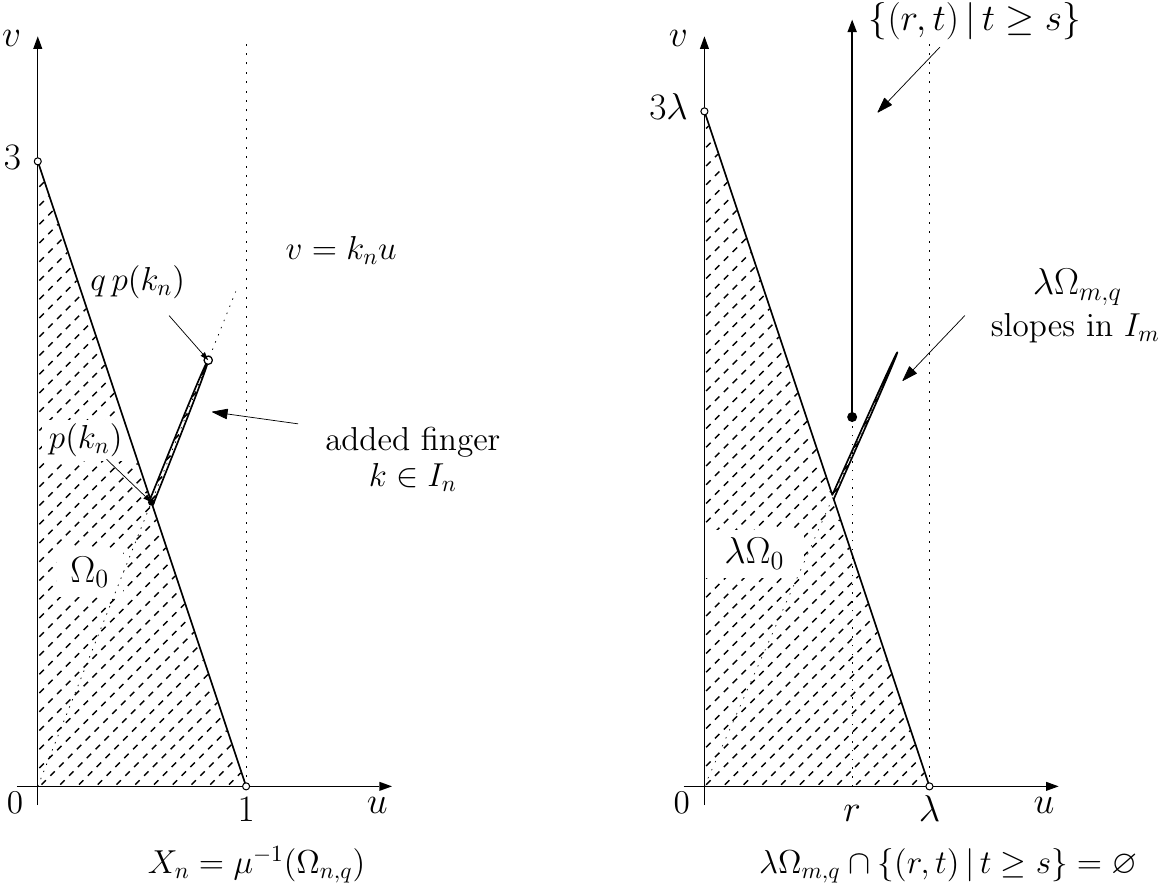}
  \caption{The construction and obstruction in Theorem~\ref{thm:G}.}
  \label{fig:theorem-G-construction-obstruction}
\end{figure}

Let us start with some preparations. For $k\in \R_{>0}$, the ray of slope $k$ meets the hypotenuse of $\mu(E(1,3))$ at $p(k):=\left(\frac3{3+k},\frac{3k}{3+k}\right)$. In the range $2<k<3$, its coordinates satisfy
\[
 \frac12<\frac3{3+k}<\frac35,\qquad
 \frac65<\frac{3k}{3+k}<\frac32,\qquad
\]
For $n\geq1$, set
\[
 k_n\coloneqq2+2^{-n},\qquad
 \delta_n\coloneqq2^{-n-3},\qquad
 I_n\coloneqq(k_n-\delta_n,k_n+\delta_n).
\]
The closures of these intervals are pairwise disjoint and contained in $(2,3)$. Moreover, $n<m$ implies $\sup I_m<\inf I_n$, so later intervals lie strictly below earlier ones in slope. Also, for each $n$, choose a smooth bump function $\chi_n: [0,\infty) \to [0,1]$ such that $\chi_n(k_n)=1$ and  vanishes outside $I_n$. Fix a scalar $1<q<5/3$,  independent of $k_n$ and $\delta_n$ (this choice will used to show that the domains $X_{n,q}\subset Z^4(1)$ in the conclusion of Lemma \ref{lem:finger-geometry} below). We define the following function of the slope $k$, 
\[
 H_{n,q}(k)\coloneqq1+(q-1)\chi_n(k).
\]
Then consider the moment region defined by
\begin{equation} \label{moment-image-stick}
 \Omega_{n,q} \coloneqq \left\{(u,v)\in\R^2_{\geq0}\;\middle|\; u >0, \,\, u+\frac v3<H_{n,q}\!\left(\frac vu\right)\right\} \cup\{(0,v)\mid 0\leq v<3\},
\end{equation}
and $X_{n,q}\coloneqq\mu^{-1}(\Omega_{n,q})$. The profile curve of $\Omega_{n,q}$ in the direction of slope $k$ is precisely $H_{n,q}(k)\cdot p(k)$. Write $\Omega_0:=\mu(E(1,3))$. Therefore, $\Omega_{n,q}$ agrees with $\Omega_0$ outside the interval of slopes $I_n$, and at slope $k_n$, it reaches $q p(k_n)$. The left picture of Figure \ref{fig:theorem-G-construction-obstruction} shows the moment image $\Omega_{n,q}$. 

\begin{lemma}\label{lem:finger-geometry}
For every $1<q<5/3$ and $n\geq1$, the toric domain $X_{n,q}$ defined via (\ref{moment-image-stick}) above is bounded and star-shaped, satisfying $E(1,3)\subset X_{n,q}\subset qE(1, 3)$, $X_{n,q}\subset Z^4(1)$, and 
\[
 \Vol(X_{n,q}) \leq\frac9{25}(q^2-1)\delta_n + \frac32.
\]
In particular, for each fixed $q$, the volume $\Vol(X_{n,q})$ approaches $3/2$ as $n \to \infty$. 
\end{lemma}

\begin{proof}
The bounds $1\leq H_{n,q}\leq q$ yield $E(1,3) \subset X_{n,q}\subset qE(1,3)$. Meanwhile, besides $E(1,3)$ inside $X_{n,q}$, every extra moment point of $X_{n,q}$ has the form $t \left(\frac3{3+k},\frac{3k}{3+k}\right)$ with $k\in I_n\subset(2,3)$ and $1\leq t<H_{n,q}(k)\leq q$. By our choice of $q$, lying inside $(1, 5/3)$, we know the first coordinate is 
\[ \frac{3t}{3+k} <\frac{3}{3+2} \cdot \frac{5}{3} = 1,\]
which confirms that $X_{n,q}\subset Z^4(1)$. Moreover, with our normalization of the moment map, symplectic volume equals the Euclidean area of the moment region. Then 
\begin{align*}
 \Vol(X_{n,q}) &= \Vol(E(1,3)) + \int_{I_n}\int_1^{H_{n,q}(k)}\frac{9t}{(3+k)^2}\,dt\,dk\\
 &=\frac{3}{2} + \frac92\int_{I_n}\frac{H_{n,q}(k)^2-1}{(3+k)^2}\,dk \leq\frac32+\frac9{25}(q^2-1)\delta_n,
\end{align*}
where we used $k>2$ and $|I_n|=2\delta_n$. This completes the proof. 
\end{proof}

\FloatBarrier

The construction yields exact distance estimates, which we now establish.

\begin{prop}\label{prop:finger-distance}
For the domains $X_{n,q}$ defined via (\ref{moment-image-stick}) above and every fixed $1<q<5/3$,
\[
 \dBM(X_{n,q},E(1,3))=\log q \,\,\,\,\mbox{and}\,\,\,\,
 \dBM(X_{n,q},X_{m,q})=\log q\quad(n\neq m).
\]
\end{prop}

\begin{proof}
The inclusions in Lemma~\ref{lem:finger-geometry} give $X_{n, q}\subset qE(1,3)\subset qX_{m, q}$ and $X_{m, q}\subset qE(1,3)\subset qX_{n, q}$. Hence both $d_{\rm BM}$ are at most $\log q$. To show the lower bound, we need to use the rigidity of Lagrangian intersection. 

Suppose $n<m$, and fix parameters $\lambda$ and $\tau$ satisfying $1\leq\lambda<\tau<q$. The point
\[
 (r,s)\coloneqq\tau\left(\frac3{3+k_n},\frac{3k_n}{3+k_n}\right)
\]
lies in $\Omega_{n,q}$, since $H_{n,q}(k_n)=q$ (and $\tau < q$). Therefore, the product Lagrangian torus $L(r,s)\subset X_{n, q}$. Moreover, by our choice $k_n = 2 + 2^{-n}$, estimating the coordinates gives
\[
 \frac\lambda2<r<1\leq\lambda<s,\qquad s=k_n r>2r.
\]
In particular, $(r,s)$ is in the region $\Cred$. Then we claim that the vertical ray $\{(r,v)\mid v\geq s\}$ is disjoint from $\lambda\Omega_{m,q}$. Indeed, it misses the rescaled triangle $\lambda \mu(E(1,3))$, because every point on the ray satisfies
\[
 r+\frac v3\geq r+\frac s3=\tau>\lambda.
\]
It also misses the extra part of $\lambda\Omega_{m,q}$, because every point in this extra part has slope in $I_m$, whereas every point on the ray has slope at least $s/r=k_n>\sup I_m$. This proves the claim.

If a symplectic embedding $X_{n, q}\into\lambda X_{m, q}$ existed, the monotonicity of the reduced shape invariant would provide a Lagrangian torus $L$ in $\lambda X_{m, q}$ with a Maslov-two basis of actions $(r,s)$. Since $\lambda X_{m, q}\subset Z^4(\lambda)$, Theorem~\ref{thm:A} applies and implies that this torus $L$ must meet
\[
 \bigcup_{v\geq s}L(r,v)=\mu^{-1}(\{(r,v)\mid v\geq s\}),
\]
contradicting the claim. So, $X_{n, q}$ does not embed into $\lambda X_{m, q}$ for any $1\leq\lambda<q$, which provides the required lower bound. A symmetric argument will cover the case when $m<n$.

Finally, for the distance to $E(1,3)$, use the same torus in $X_{n, q}$. The ray $\{(r,v)\mid v\geq s\}$ misses $\lambda \mu(E(1,3))$, and $\lambda E(1,3)\subset Z^4(\lambda)$. The same application of Theorem~\ref{thm:A} prohibits the embedding $X_{n, q}\into\lambda E(1,3)$ for every $1\leq\lambda<q$, proving that $\dBM(X_{n, q},E(1,3))\geq\log q$.
\end{proof}

Now we are ready to prove Theorem~\ref{thm:G}. 

\begin{proof}[Proof of Theorem~\ref{thm:G}]
Given $0<\ep<\log(5/3)$, set $q\coloneqq e^\ep$, then the toric domains $X_n:= X_{n,q}$ defined via (\ref{moment-image-stick}) provide the desired family by Lemma~\ref{lem:finger-geometry} and Proposition~\ref{prop:finger-distance} (in particular, their volumes converge to $3/2$).

Now, let us make modifications so that our domains have the same fixed volume. The scaling convention \eqref{eq:scaling} and volume monotonicity give, for all $U \in \mathcal S$ and $a>0$,
\[
 \dBM(U,aU)=|\log a|.
\]
Then choose $0<\ep<\log(5/3)$, and take the sequence $X_n = X_{n, q}$ as above corresponding to $q=e^\ep$. Define
\[
 \alpha_n\coloneqq\left(\frac{3/2}{\Vol(X_n)}\right)^{1/2} \,\,\,\,\mbox{and}\,\,\,\,
 Y_n\coloneqq\alpha_n X_n.
\]
Then $\Vol(Y_n)=3/2$, $0<\alpha_n\leq1$ and the scalars $\alpha_n$ converge to $1$, and all $Y_n$ remain star-shaped. Moreover, the triangle inequality yields
\[
\begin{aligned}
 \dBM(Y_n,E(1,3))&\leq\ep+|\log\alpha_n|=\ep-\log\alpha_n,\\
 \dBM(Y_n,Y_m)&\geq\ep-|\log\alpha_n|-|\log\alpha_m|\\
 &=\ep+\log\alpha_n+\log\alpha_m\qquad(n\neq m).
\end{aligned}
\]
By the volume bound in Lemma \ref{lem:finger-geometry}, after discarding finitely many terms we may assume $|\log\alpha_n|<\ep/4$. Thus the remaining $Y_n$ lie within $d_{\rm BM}$-distance $5\ep/4$ from $E(1,3)$ and have pairwise $d_{\rm BM}$-distance greater than $\ep/2$. This proves the local non-compactness at $E(1,3)$ (in the volume-$3/2$ subspace of $\mathcal S$).\end{proof}

\

\bibliographystyle{amsplain}
\bibliography{Sh_ref}

\

\end{document}